\documentclass[11pt,a4paper]{article}
\usepackage[utf8]{inputenc}
\usepackage{amsmath}
\usepackage{amsfonts}
\usepackage{amssymb}
\usepackage{tikz}
\usepackage{graphicx}
\usepackage{amsthm}
\usepackage[english]{babel}
\usepackage{array}

\usepackage{caption}

\usetikzlibrary{matrix}
\usepackage{enumitem}

\theoremstyle{plain}
\newtheorem{The}{Theorem}[section]
\newtheorem*{theorem*}{Theorem}

\newtheorem{Lem}[The]{Lemma}
\newtheorem*{Definition}{Definition}
\newtheorem{Prop}[The]{Proposition}

\newtheorem{Cor}[The]{Corollary}

\newcommand{\bl}{\begin{Lem}}
\newcommand{\el}{\end{Lem}}

\newcommand{\bd}{\begin{Definition}}
\newcommand{\ed}{\end{Definition}}

\newcommand{\bt}{\begin{The}}
\newcommand{\et}{\end{The}}

\newcommand{\bp}{\begin{Prop}}
\newcommand{\ep}{\end{Prop}}

\newcommand{\bc}{\begin{Cor}}
\newcommand{\ec}{\end{Cor}}

\newcommand{\g}{\mathfrak{g}}

\newcommand{\kk}{\mathfrak{k}}

\newcommand{\p}{\mathfrak{p}}

\newcommand{\h}{\mathfrak{h}}

\newcommand{\hp}{\h_{\p}}
\newcommand{\Ad}{\text{Ad}}
\newcommand{\ad}{\text{ad}}
\newcommand{\tx}{\mathbf{x}}
\newcommand{\G}{G\slash K}
\newcommand{\tG}{\tilde{G}\slash\tilde{K}}
\newcommand{\tM}{\tilde{M}}

\newcommand{\R}{\mathbb{R}}
\newcommand{\C}{\mathbb{C}}
\newcommand{\AU}{\mathrm{A}(eK)}
\newcommand{\AT}{\mathrm{A}_T(eK)}
\newcommand{\Z}{\mathbb{Z}}

\newcommand{\spa}{\hspace{.2cm}}
\newcommand{\vsk}{\vskip .2cm}

\begin{document}
\title{A complete description of the antipodal set of most symmetric spaces of compact type}
\author{Jonas Beyrer}

\maketitle

\textbf{Abstract:} It is known that the antipodal set of a Riemannian symmetric space of compact type $\G$ consists of a union of $K$-orbits. We determine the dimensions of these $K$-orbits of most irreducible symmetric spaces of compact type. The symmetric spaces we are not going to deal with are those with restricted root system $\mathfrak{a}_r$ and a non-trivial fundamental group, which is not isomorphic to $\Z_2$ or $\Z_{r+1}$. For example, we show that the antipodal sets of the Lie groups $Spin(2r+1)\:\: r\geq 5$, $E_8$ and $G_2$ consist only of one orbit which is of dimension $2r$, 128 and 6, respectively; $SO(2r+1)$ has also an antipodal set of dimension $2r$; and the Grassmannian $Gr_{r,r+q}(\mathbb{R})$ has a $rq$-dimensional orbit as antipodal set if $r\geq 5$ and $q>0$.\\

\let\thefootnote\relax\footnote{\textbf{M.S.C. 2010:} 22E46, 53C30.}
Key words: antipodal set, Riemannian symmetric spaces, cut locus\\

Acknowledgment: I want to thank Oliver Goertsches very much for his many helpful comments.

\section{Introduction}

The \emph{antipodal set} of a point $p\in M$ in a connected, compact Riemannian manifold $(M,g)$ is the set of points $q\in M$ with maximal distance to $p$ and is denoted by $\mathrm{A}(p)$. Given a general Riemannian manifold $M$, it is not known how to determine $\mathrm{A}(p)$, but if $M$ is a symmetric space of compact type the situation changes; bringing us to the purpose of this paper: We give a complete description of the antipodal set of most (Riemannian) symmetric space of compact type.

The problem of determining the antipodal set in a symmetric space is not new. Already in 1978 J. Tirao \cite{Ti} solved it for symmetric spaces of rank one. However for higher rank symmetric spaces little was known. Only S. Deng and X. Liu \cite{DL} determined those compact simply connected symmetric spaces that have a finite number of points as antipodal set. In addition they gave the exact number of points.

We note that, given a symmetric space of compact type, it is enough to determine the antipodal set in each of the irreducible components, as by de Rham's decomposition theorem the antipodal set of the whole space is the product of the antipodal sets in the irreducible components. Therefore we assume from now on $M\cong\G$ to be an irreducible symmetric space of compact type, if not stated otherwise. Being able to restrict to irreducible spaces is important for us, as we do case by case calculations. 

It is quite immediate that the antipodal set of $p\in \G$ is a union of orbits of the form $\exp_{p}(\Ad(K)\textbf{x})$ for specific $\textbf{x}\in\p\cong T_p \G$.
L. Yang determined in \cite{Yang} those $\textbf{x}$ of all irreducible compact simply connected symmetric spaces and in \cite{Yang2} those $\textbf{x}$ of most irreducible non-simply connected symmetric spaces of compact type. He did not determine $\tx$ for those $\G$ that have restricted root system $\mathfrak{a}_r$ and a non-trivial fundamental group, which is not isomorphic to $\Z_2$ or $\Z_{r+1}$. Since the knowledge of $\tx$ is essential for this work, we are not able to determine the antipodal set in those cases. Therefore we refer to them as \emph{excluded} cases. 

In this paper we analyze the antipodal set of $\G$ by determining the dimensions of the orbits $\exp_{eK}(\Ad(K)\textbf{x})$ building the antipodal set, if $\G$ is not one of the excluded cases.

The paper is organized as follows: In section 2 we describe some well known facts on symmetric spaces and their antipodal set. In section 3 we use those facts to give an explicit description of the tangent space of the antipodal set of $\G$. If $\G$ is in addition simply connected this description reads as follows: Let $\alpha_i$ be simple roots of the restricted root system $\Sigma$ of $\G$, $\psi=\sum_{i=1}^r d_i \alpha_i$ the highest root, $J_j^{'}=\lbrace \alpha =\sum_{k=1}^{r} c_k \alpha_k \in \Sigma^{+}\:|\: \frac{c_j}{d_j}\notin \mathbb{N} \rbrace$, $\p(\alpha)$ the root spaces in $\p$ and $\textbf{x}=\pi e_j$, where $e_j$ is a maximal corner. Then the tangent space of $\exp_{eK}(\Ad(K)\textbf{x})\subset \mathrm{A}(eK)$ is a parallel translate of $\bigoplus_{\alpha \in J_j^{'}} \p(\alpha)$. For a precise description of maximal corners see section \ref{Yang-section}.
In section 4 we use this description to determine the dimensions of all orbits in all $\G$ but the excluded ones explicitly case be case. We give several example calculations, the other cases can be treated in the same manner. In the tables at the end of section \ref{section-dimensions} the whole results are stated. 

We want to remark that the natural numbers $\mathbb{N}$ contain zero in our notation and furthermore that a different definition of antipodal set of symmetric spaces exists in the literature, describing a different object.

\section{Preliminaries}

\subsection{Some facts on symmetric spaces of compact type}

Before we start with analyzing the problem, we want to coarsely remind of some properties of symmetric spaces of compact type. A proper description of the following and definitions of the named objects can for example be found in \cite{Hel}.\\

A compact simply connected symmetric spaces is of compact type. Furthermore the simply connected cover of a symmetric space of compact type is compact itself.

An irreducible symmetric space of compact type $M$ is isometric to the space $G\slash K:=I_0 (M)\slash (I_0(M))_p$ with $I_0 (M)$ being the connected component of the identity of the isometry group of $M$ and $(I_0(M))_p$ being the stabilizer of a point $p\in M$ of the natural action of the isometry group. 

The metric on $G\slash K$ corresponds to a left invariant extension of a multiple of the Killing form on $Lie(G)=\g$. As $\G$ is of compact type, $\g$ is semi-simple. Furthermore there exists a natural involution on $\g$. Let $\mathfrak{k}$ and $\mathfrak{p}$ be the 1 and -1 eigenspaces of this involution, respectively. The decomposition $\g=\mathfrak{k}+\mathfrak{p}$ is called the \emph{Cartan decomposition}. The space $\mathfrak{p}$ has the property that it is isometric to $T_{eK} G\slash K$. 

Let $\hp$ be a maximal abelian subspace in $\p$. We define for $\alpha\in \hp^{*}$
\begin{align*}
\g(\alpha)&:=\lbrace X\in\g\:|\:[H,[H,X]]=-\alpha(H)^2 X\hspace{.2cm} \text{for all }H\in \hp\rbrace.
\end{align*}
If $\g(\alpha)\neq 0$ and $\alpha\neq 0$, we call $\alpha$ a \emph{root}. Furthermore we set $\kk(\alpha):=\g(\alpha)\cap\kk$, $\p(\alpha):=\g(\alpha)\cap\p$ . We call $\dim \kk(\alpha)$ the \emph{multiplicity} of the root $\alpha$ and the set of all roots $\Sigma$ \emph{restricted root system}.

Let $B(\cdot,\cdot)$ be the Killing from on $\g$. As $\g$ is semi-simple, we can identify an element $\alpha\in\hp^{*}$ with an element $H_{\alpha}\in\hp$ by the relation $\alpha(H)=B(H_{\alpha},H)$ for all $H\in\hp$. We define $B$ on $\hp^{*}$ by $B(\alpha,\beta):=B(H_{\alpha},H_{\beta})$ for $\alpha,\beta\in\hp^{*}$. In this case $(\Sigma, c B(\cdot,\cdot))$ is an abstract root system for an suitable constant $c\in\R$. For notational reasons we just write $(\cdot,\cdot)$ instead of $cB(\cdot,\cdot)$. For the abstract root system we can fix an ordering and get a set of positive roots $\Sigma^{+}$. With respect to $\Sigma^{+}$ we can decompose the spaces $\kk$ and $\p$ as follows
\begin{align}\label{decompositionKandP}
\kk=\kk(0)\oplus\bigoplus\limits_{\alpha\in\Sigma^{+}}\kk(\alpha),\hspace{1cm}
\p=\p(0)\oplus\bigoplus\limits_{\alpha\in\Sigma^{+}}\p(\alpha).
\end{align}
These decompositions are called \emph{root space decomposition} of $\kk$ and $\p$, respectively. Let $\Sigma$ be a root system, a set $\tilde{\Sigma}\subset\Sigma$ is called \emph{root subsystem}, if the following holds:

\begin{enumerate}
\item for $\alpha,\beta\in\tilde{\Sigma}$ it is $\alpha +\beta\in\tilde{\Sigma}$, if $\alpha +\beta\in\Sigma$,
\item $-\tilde{\Sigma}=\tilde{\Sigma}$.
\end{enumerate}

A root subsystem is an abstract root system itself and clearly there is a unique ordering compatible with the ordering on $\Sigma$. Similar to the root space decompositions we define for a root subsystem $\tilde{\Sigma}$
\begin{align}\label{decomposition-subsystem}
\kk(\tilde{\Sigma}):=\kk(0)\oplus \bigoplus\limits_{\alpha\in\tilde{\Sigma}^{+}} \kk(\alpha).
\end{align}

\subsection{Basic properties of the antipodal set}

Mathematically, the antipodal set of a point $p\in G\slash K$ is given by
\begin{align*}
\mathrm{A}(p) =\lbrace x\in G\slash K \:|\: d(x,p)\geq d(y,p)\hspace{.2cm}  \forall y\in G\slash K \rbrace.
\end{align*} 
Let $p=aK\in G\slash K$ be arbitrary and let $\ell_a$ be the left multiplication by $a$. Then
\begin{align*}
d(y,eK)\leq d(x,eK)\hspace{.2cm} \hspace{.2cm} \Longleftrightarrow\hspace{.2cm} d(\ell_a(y),p)\leq d(\ell_a(x),p)\hspace{.2cm}  \forall x,y\in G\slash K,
\end{align*}
which implies $\mathrm{A}(p)=\ell_a(\AU )$. Let $d(\G)$ denote the diameter. As $\G$ is compact, we find $p_0,q_0\in \G$ with $d(p_0,q_0)=d(\G)$. The isometry group of a symmetric space acts transitively therefore $q\in A(p)$ if and only if $d(p,q)=d(\G)$. Hence every $q\in \mathrm{A}(p)$ can be joined to $p$ by a geodesic which is minimizing till $q$ but not beyond. It follows $\mathrm{A}(p)\subset C(p)$, where $C(p)$ denotes the cut locus of $p$. Thus $\mathrm{A}(p) =\lbrace x\in C(p) \:|\:d(x,p)\geq d(y,p)\:\:\forall y\in C(p)\rbrace$.\\
The set of points $X\in T_p \G$ with $\exp_p (tX)$ being a minimizing geodesic for $t<1$ and not minimizing for $t>1$ is called the \emph{cut locus of $p$ in $T_p \G$} or alternatively \emph{the cut locus of $p$ in the tangent space} and is denoted by $C_{T} (p)$. We use this to define the \emph{antipodal set in the tangent space} to be
\begin{align}\label{antipodalset-tangetspace}
\mathrm{A}_T(p) :=\lbrace X\in C_T(p) \:|\: |X|\geq |Y| \:\:\forall\: Y\in C_T(p) \rbrace.
\end{align}
If $X\in C_T(p)$, then $|X|=d (\exp_{p}(X),p)$ and $\exp_p$ maps $C_T(p)$ onto $C(p)$, thus
\begin{align}\label{antipodal-tangent-regular}
\exp_{p}(\mathrm{A}_T(p))= \mathrm{A}(p).
\end{align}
As the antipodal sets of two points are isometric, it is enough to consider only $eK$. In view of  
\eqref{antipodal-tangent-regular} we need to make two steps, namely determining $\mathrm{A}_T(eK)$ and analyzing $\exp_{eK}$.

\subsection{Cut locus and antipodal set in the tangent space for simply connected symmetric spaces}\label{Yang-section}

Let $\Sigma$ be the restricted root system of $G\slash K$, where $\G$ is an irreducible compact simply connected symmetric space. We fix an ordering on $\Sigma$ that we keep for the rest of this paper. Furthermore let $\Sigma^{+}$ be the set of positive roots, $\Pi^{\Sigma}$ the set of simple roots and $\psi$ the highest roots of $\Sigma$ with respect to this ordering. There is only one highest root, as $\G$ is irreducible. L. Yang \cite{Yang} defined the \emph{Cartan polyhedron} to be
\begin{align}\label{Cartan-polyhedron}
\triangle:=\lbrace x\in \h_{\p}\:|\: \gamma(x)\geq 0 \text{ for } \gamma \in \Pi^{\Sigma}\:\wedge\: \psi(x)\leq 1\rbrace.
\end{align} 
Clearly the condition that $\gamma(x)\geq 0$ for $\gamma \in \Pi^{\Sigma}$ implies that $\triangle$ is contained in the closure of a Weyl chamber. The side of the Cartan polyhedron that does not contain 0 is denoted by\vskip -0.3cm

\begin{align}\label{cartanpolyhedronedge}
\triangle^{'}:=\lbrace x\in \triangle\:|\: \psi(x)=1\rbrace.
\end{align}

Now we are able to cite the following theorem, which gives an explicit description of the cut locus in the tangent space.

\bt\label{CartanPMain}\label{Yang-theo}(see \cite{Yang} p. 689 and Appendix)\\
Let $G\slash K$ be a compact simply connected symmetric space. Then $C_T(eK)=\Ad(K)(\pi \triangle^{'})$.
\et

An immediate consequence of this theorem is, by \eqref{antipodalset-tangetspace} and the fact that $K$ acts orthogonally, that\vskip -.7cm
\begin{align*} 
\AT=\Ad(K)(\pi \:\text{max}(\triangle^{'})),
\end{align*}
with $\max(\triangle^{'}):=\lbrace x\in \triangle^{'}\:|\: |x|\geq |y| \:\forall y\in \triangle^{'}\rbrace$. 

The following proposition is well known.

\bp\label{K orbit Weyl chamber intersection}\mbox{}\\
Let $\G$ be an irreducible symmetric space of compact type. For an orbit of the form $\Ad(K)\tx$, with $\tx\in \bar{C}$ and $\bar{C}$ being a closed Weyl chamber, we have $\Ad(K)\tx\cap \bar{C}=\lbrace\tx\rbrace$.
\ep

This proposition implies that $\AT$ consists of a union of orbits and the set of orbits is in one to one correspondence to the set $\text{max}(\triangle^{'})$. L. Yang \cite{Yang} has determined the set $\max(\triangle^{'})$ more precise. For the convenience of the reader, we explain his approach and some of his steps:\\ 
As said before, $(\Sigma,(\cdot,\cdot))$ is an abstract irreducible root system. Hence the Weyl chamber is a cone and we only have one highest root. It follows that the Cartan polyhedron is a simplex. We observe that for $x_1,x_2\in \triangle$ and $t\in [0,1]$, we have
\begin{align}\label{maximal elements are in the croners}
 &(tx_1+(1-t)x_2,tx_1+(1-t)x_2)^{\frac{1}{2}}\leq \max \lbrace (x_1,x_1)^{\frac{1}{2}},(x_2,x_2)^{\frac{1}{2}}\rbrace.
\end{align}

We conclude that the function $\triangle\to\R$ sending $x\mapsto (x,x)^{\frac{1}{2}}$ reaches its maximum on the vertices, but not on $0$. Let $\lbrace \alpha_1,\ldots,\alpha_r\rbrace$ be a set of simple roots of $\Sigma$ and $\psi=\sum_{i=1}^r d_i \alpha_i$ its highest root. Each side of the Cartan polyhedron is contained in a hyperplane of the form $\lbrace x\in \hp \:|\: \alpha_i(x)=0\rbrace$ or $\lbrace x\in \hp\:|\: \psi(x)=1\rbrace$, hence the vertices of the Cartan polyhedron are the points that lie in $r$ of these $r+1$ hyperplanes. These points  are $0$ and $e_j$ with $j=1,\ldots,r$ such that\vskip -.6cm
\begin{align}\label{Corners-polyh-equation}
\alpha_i(e_j)=\frac{1}{d_j}\delta_{ij}.
\end{align}\vskip -.2cm
In particular we get
\begin{align*}
\max (\triangle^{'})=\lbrace e_j \:|\: (e_j,e_j)^{\frac{1}{2}}=\max_{i=1,\ldots,r} (e_i,e_i)^{\frac{1}{2}}\rbrace.
\end{align*} 
We call the corners $e_j$ of the Cartan polyhedron with $e_j\in\max (\triangle^{'})$ \emph{the maximal corners of the Cartan polyhedron}. For each irreducible restricted root system those corners can be determined with straight forward calculations. This is done in \cite{Yang}. The results are listed in the following table. We want to remark that we used the indexing of the simple roots as in \cite{Hel} p. 477, 478, which differs in the cases $\mathfrak{e}_6, \mathfrak{e}_7, \mathfrak{e}_8$ and $\mathfrak{g}_2$ from the one in \cite{Yang}. The factors $d_j$ of the highest root are well known (see \cite{Bou2} Plate I - IX).\vskip .3cm

\captionof{table}{Maximal corners of the Cartan polyhedron and the corresponding factors of the highest root (see \cite{Yang} p. 689 - 693)} \label{table-maximal-corners-and-factors} 
\begin{tabular}{|p{1.5cm}|p{1.4cm}|p{2.85cm}|p{1.5cm}|p{1.4cm}|p{2.7cm}|}
\hline
$\Sigma$ & $\max(\triangle^{'})$ & \begin{footnotesize} Factors \end{footnotesize} $d_j$ & $\Sigma$ & $\max(\triangle^{'})$ & \begin{footnotesize} Factors \end{footnotesize} $d_j$   \\ \hline
$\mathfrak{a}_{2r}$ & $e_r; e_{r+1}$  & $d_r=1; d_{r+1}=1$ & $\mathfrak{d}_r, r>4$ & $e_{r-1}$; $e_r$ & $d_{r-1}=1; d_r=1$\\ \hline
$\mathfrak{a}_{2r-1}$ & $e_r$ & $d_r=1$ & $\mathfrak{e}_6$ & $e_1; e_6$ & $d_1=1; d_6=1$ \\ \hline
$\mathfrak{b}_2, \mathfrak{b}_3$ & $e_1$ & $d_1=1$ & $\mathfrak{e}_7$ & $e_7$ & $d_7=1$\\ \hline
$\mathfrak{b}_4$ & $e_1; e_4$ & $d_1=1; d_4=2$ & $\mathfrak{e}_8$ & $e_1$ & $d_1=2$ \\ \hline
$\mathfrak{b}_r, r>4$ & $e_r$ & $d_r=2$ & $\mathfrak{f}_2$ & $e_4$ & $d_4=2$\\ \hline
$\mathfrak{c}_r$ & $e_r$ & $d_r=1$ & $\mathfrak{g}_2$ & $e_1$ & $d_1=3$\\ \hline 
$\mathfrak{d}_4$ & $e_1; e_3; e_4$ & $d_1=d_3=d_4=1$ & $\mathfrak{bc}_r$ & $e_r$ & $d_r=2$ \\ \hline
\end{tabular}
The maximal corners are not explicitly stated in \cite{Yang}, but can easily be deduced.\vskip .4cm

Exemplary, we give L. Yangs calculation for $\Sigma=\mathfrak{a}_r$: A choice of simple roots of $\Sigma$ is $\alpha_1=x_1-x_2,\ldots,\alpha_r=x_r-x_{r+1}$ for $\lbrace x_1,x_2,\ldots,x_{r+1}\rbrace$ a basis of $\R^{r+1}$ such that $(x_i,x_j)=\frac{1}{2} (\psi,\psi)\delta_{ij}$ and $\psi=\sum_{i=1}^r \alpha_i$ the highest root (see \cite{Bou2} p. 265). We see that the factors $d_i$ in front of the simple roots building $\psi$ are all equal to 1. As the corners of the Cartan polyhedron satisfy $\alpha_i(e_j)=(\alpha_i,e_j)=\frac{1}{d_j}\delta_{ij}$, we can deduce
\begin{align*}
e_j=\frac{2}{(\psi,\psi)(r+1)}((r+1-j)\sum\limits_{k=1}^j x_k-j\sum\limits_{k=j+1}^{r+1} x_k)\hspace{.7cm} 1\leq j \leq n.
\end{align*}
This gives
\begin{align*}
(e_j,e_j)=\frac{2j(r+1-j)}{(\psi,\psi)(r+1)}.
\end{align*}
For $r$ odd this is maximal if $j=\frac{r+1}{2}$ and for $r$ even this is maximal if $j=\frac{r}{2}$ or $j=\frac{r}{2}+1$. Hence the maximal corner of the Cartan polyhedron corresponding to $\mathfrak{a}_{2r-1}$ is $e_r$ and for $\mathfrak{a}_{2r}$ the maximal corners are $e_r, e_{r+1}$.

\subsection{Cut locus and antipodal set in the tangent space for non-simply connected symmetric spaces}

Let $\tilde{M}$ be an irreducible compact simply connected symmetric space. We can write $\tilde{M}=\tilde{G}\slash \tilde{K}$ with $\tilde{G}$ simply connected. Throughout this section $\tilde{G}$ is always chosen simply connected. We define 
\begin{align*} 
Z_{\tG}(\tilde{K}):=\lbrace p\in \tG \:|\: k\cdot p=p \:\:\:\forall k\in \tilde{K}\rbrace. 
\end{align*}
For this set L. Yang showed the following:

\bp (see \cite{Yang2} p. 517 - 519)\\
Notations as before. For every $p=a\tilde{K}\in \tilde{G}\slash\tilde{K}\backslash \lbrace e\tilde{K}\rbrace$ the following are equivalent:
\begin{enumerate}[label=(\alph*)]
\item $p\in Z_{\tG}(\tilde{K})$;
\item $a\in N_{\tilde{G}}(\tilde{K})$, where $N_{\tilde{G}}(\tilde{K})$ denotes the normalizer of $\tilde{K}$ in $\tilde{G}$;
\item $p= \exp_{e\tilde{K}} (\pi e_j)$, where $e_j$ is a corner of the Cartan polyhedron such that $d_j =1$.
\end{enumerate}
Furthermore $Z_{\tG}(\tilde{K})$ is a finite abelian group which can be identified with a subgroup of $\tilde{G}$.
\ep 

We want to give some short comments on this proposition. The equivalence of $(a)$ and $(b)$ can be easily seen from the fact that $ka\tilde{K}=a\tilde{K}$ for all $k\in\tilde{K}$ and $a\tilde{K}\in Z_{\tG}(\tilde{K})$. From this it follows that $Z_{\tG}(\tilde{K})$ is a group. Part of the proof is to show that the map $\Psi(a\tilde{K}):=a \sigma^{-1}(a)$ is well defined on $Z_{\tG}(\tilde{K})$, where $\sigma$ is the natural involution on $\tilde{G}$. Furthermore one shows that $\Psi: Z_{\tG}(\tilde{K})\to Z(\tilde{G})$ is a monomorphism, while $Z(\tilde{G})$ denotes the center of $\tilde{G}$. Then $Z_{\tG}(\tilde{K})$ is a finite abelian group and the map $\Psi$ gives an identification with a subgroup of $\tilde{G}$.

\bp (see \cite{Yang2} p. 519)\\
Every symmetric space $M$ of compact type with simply connected cover $\tM=\tG$ can be expressed as a quotient $\tM\slash \Gamma$ for a subgroup $\Gamma < Z_{\tM}(\tilde{K})$ and for every subgroup $\Gamma < Z_{\tM}(\tilde{K})$ the space $\tM\slash \Gamma$ is a symmetric space which is covered by $\tM$.
\ep 

The quotient $M=\tM\slash\Gamma$ is called a \emph{Clifford-Klein form of $M$}.\\
Let $\Sigma$ be the restricted root system of $\tM$ and hence also of $M$. In addition let $\psi$ be the highest root, $\triangle$ the Cartan polyhedron and $e_i$ the corners of the Cartan polyhedron. For a symmetric space $M$ with Clifford-Klein-form $\tM\slash\Gamma$ we define, same as L. Yang \cite{Yang2}, the sets 
\begin{align*}
P_{\Gamma}&:=\lbrace x\in \triangle\:|\: (x,e_i)\leq \frac{1}{2}(e_i,e_i) \text{ for every } \exp_{e\tilde{K}}(\pi e_i)\in\Gamma \rbrace,\\
P_{\Gamma}^{'}&:=\lbrace x\in P_{\Gamma}\:|\: (x,\psi)=1\: \vee\: \exists \exp_{e\tilde{K}}(\pi e_i)\in \Gamma: (x,e_i)=\frac{1}{2}(e_i,e_i)\rbrace.
\end{align*} 

The set $P_{\Gamma}$ or more precisely its boundaries which do not contain zero, namely $P_{\Gamma}^{'}$, play an important role in the description of the cut locus in the tangent space. This is given by the next theorem.

\bt\label{Yang2-theorem} (see \cite{Yang2} p. 521 and Appendix)\\
Let $\tM=\tG$ be an irreducible compact simply connected symmetric space and let $M=\G$ be a symmetric spaces covered by $\tM$. Let $M=\tM\slash \Gamma$ be a Clifford-Klein form, where $\Gamma$ is a non-trivial subgroup of $Z_{\tM}(\tilde{K})$. Then 
\begin{align*}
C_T(eK)=\Ad(K)(\pi P_{\Gamma}^{'}).
\end{align*}
\et 

Similar as in the simply connected case, this theorem implies together with \eqref{antipodalset-tangetspace} and the fact that $K$ acts orthogonally that\vskip -.7cm
\begin{align*} 
\AT=\Ad(K)(\pi \:\text{max}(P_{\Gamma}^{'})),
\end{align*}
where $\max(P_{\Gamma}^{'}):=\lbrace x\in P_{\Gamma}^{'}\:|\: |x|\geq |y| \hspace{.2cm}\forall y\in P_{\Gamma}^{'}\rbrace$ and $\G$ is non-simply connected. By Proposition \ref{K orbit Weyl chamber intersection} we can derive a one to one correspondence of points in  $\max(P_{\Gamma}^{'})$ and orbits building the antipodal set in the tangent space.\\
We want to determine the set $\max(P_{\Gamma}^{'})$. By a similar argument as in \eqref{maximal elements are in the croners} it follows that this set is a subset of the corners of the polyhedron $P_{\Gamma}$. The explicit description of $\Gamma$ given at the beginning of this section, allows to try to determine $\max(P_{\Gamma}^{'})$ of every irreducible non-simply connected symmetric space in case by case calculations. L. Yang has done this implicitly in \cite{Yang2}, where he determined the diameter of irreducible non-simply connected symmetric spaces of compact type. As part of this calculations he also determined all possible subgroups of $Z_{\tM}(\tilde{K})$. In the case that $\Sigma=\mathfrak{a}_r$ and $\Gamma$ is not isomorphic to $\Z_2$ or $\Z_{r+1}$ he did not determine the diameter and therefore we don't know the set $\max(P_{\Gamma}^{'})$. However, for all the other cases we can derive $\max(P_{\Gamma}^{'})$. The results are given in the table below.
 
In the table the subgroup $\Gamma$ is given only up to isomorphism, if none of the subgroups of $Z_{\tM}(\tilde{K})$ are isomorphic. Otherwise they are given explicitly. Our indexing of the roots is again as in \cite{Hel} p. 477, 478 and differs therefore for $\mathfrak{e}_6$ and $\mathfrak{e}_7$ from the one in \cite{Yang2}. Again the factors $d_j$ of the highest root are well known.
\newpage

\captionof{table}{Maximal corners of $P_{\Gamma}$ for most non-simply connected $\G=\tilde{M}\slash\Gamma$ and the corresponding factors of the highest root (see \cite{Yang2} p. 528 - 533 and Appendix)} \label{table-non-cs-max corners} 
\setlength{\extrarowheight}{2.4pt}
\begin{tabular}{|p{3.8cm}|p{1.5cm}|p{3cm}|p{2.5cm}|} 
\hline
$\Sigma$ & $\Gamma$ & $\max(P_{\Gamma}^{'})$ & Factors of $\psi$\\ \hline  
$\mathfrak{a}_r$ \hspace{.3cm} $r\geq 3$, $r$ odd and & $\mathbb{Z}_2$ & $e_{\frac{r+1}{4}}$ & 1 \\
\hspace{.75cm} $\frac{r+1}{2}$  even &  & &  \\
$\mathfrak{a}_r$ \hspace{.3cm} $r\geq 3$, $r$ odd and & $\mathbb{Z}_2$ & $\frac{1}{2}(e_{\frac{r-1}{4}}+e_{\frac{r+3}{4}})$ & $(1,1)$ \\
\hspace{.75cm} $\frac{r+1}{2}$  odd &  & &  \\ \hline
$\mathfrak{a}_r$ & $\mathbb{Z}_{r+1}$ & $\frac{1}{r+1}(e_1+\ldots +e_r)$ & $(1,\ldots,1)$ \\ \hline
$\mathfrak{a}_r$ & otherwise & unknown &  \\ \hline
$\mathfrak{b}_r$  & $\mathbb{Z}_2$ &$e_r$ & 2\\ \hline 
$\mathfrak{c}_r$  \hspace{.3cm} $r$ even & $\mathbb{Z}_2$ & $e_{\frac{r}{2}}$ & 2 \\
$\mathfrak{c}_r$  \hspace{.3cm} $r$ odd & $\mathbb{Z}_2$ & $\frac{1}{2}(e_{\frac{r-1}{2}}+e_{\frac{r+1}{2}})$ & $(2,2)$ \\ \hline
$\mathfrak{d}_r$  \hspace{.3cm} $r$ even & $\mathbb{Z}_2\oplus \mathbb{Z}_2$ & $e_{\frac{r}{2}}$ & 2 \\
$\mathfrak{d}_r$  \hspace{.3cm} $r$ odd & $\mathbb{Z}_4$ & $\frac{1}{2}(e_{\frac{r-1}{2}}+e_{\frac{r+1}{2}})$ & $(2,2)$ \\ \hline
$\mathfrak{d}_r$  \hspace{.3cm} $r$ & $\lbrace e,p_1\rbrace$ & $e_{r-1}$; $e_r$ & 1; 1 \\ \hline
$\mathfrak{d}_r$  \hspace{.3cm} $r$ even, $r\leq 6$ & $\lbrace e,p_{r-1}\rbrace$ & $e_1$ & 1 \\
$\mathfrak{d}_8$  \hspace{.3cm} & $\lbrace e,p_{r-1}\rbrace$ & $e_1; e_4$ & 1; 2 \\
$\mathfrak{d}_r$  \hspace{.3cm} $r$ even, $r\geq 10$ & $\lbrace e,p_{r-1}\rbrace$ & $e_{\frac{r}{2}}$ & 2 \\ \hline
$\mathfrak{d}_r$  \hspace{.3cm} $r$ even, $r\leq 6$ & $\lbrace e,p_{r}\rbrace$ & $e_1$ & 1 \\
$\mathfrak{d}_8$  \hspace{.3cm} & $\lbrace e,p_{r}\rbrace$ & $e_1; e_4$ & 1; 2 \\
$\mathfrak{d}_r$  \hspace{.3cm} $r$ even, $r\geq 10$ & $\lbrace e,p_{r}\rbrace$ & $e_{\frac{r}{2}}$ & 2 \\ \hline
$\mathfrak{e}_6$ & $\mathbb{Z}_3$ & $e_4$ & 3\\ \hline
$\mathfrak{e}_7$ & $\mathbb{Z}_2$ & $e_2$ & 2 \\ \hline
\end{tabular}\vskip .1cm
We use the notation $e=e\tilde{K}$, $p_i:= \exp_{e\tilde{K}}(\pi e_i)$. In the last column we write $(d_i,d_j)$, if the maximal corner is of the form $c(e_i + e_j)$ for some $c\in\R$.\\

We want to remark that if the restricted root system $\Sigma$ is one of $\mathfrak{bc}_r, \mathfrak{e}_8, \mathfrak{f}_4$ or $\mathfrak{g}_2$, then $Z_{\tM}(\tilde{K})=\lbrace e\tilde{K}\rbrace$ and hence there is no non-simply connected symmetric space of compact type with one of those restricted root systems.\\

We give an example of how Yang rather implicitly determines $\max(P_{\Gamma}^{'})$ in the case that $\Sigma=\mathfrak{d}_r$ and $\Gamma=\lbrace e\tilde{K},\exp_{e\tilde{K}}(\pi e_1)\rbrace$: At first he gives explicit descriptions of the maximal corners of the Cartan polyhedron $e_i$ up to a scaling by $\frac{1}{2}(\psi,\psi)$, where $\psi$ is the highest root of $\mathfrak{d}_r$. From that he shows $(e_1,e_1)=2(\psi,\psi)^{-1}$ and $(e_1,e_i)=(\psi,\psi)^{-1}$ for all  $2\leq i\leq r$, but then $(e_1,e_i)=\frac{1}{2}(e_1,e_1)$. This implies that the corners of $P_{\Gamma}$ include $e_2,\ldots e_r$. By definition $P_{\Gamma}\subset \triangle$. As now $\max(\triangle^{'})=\lbrace e_{r-1},e_r\rbrace\subset P_{\Gamma}$, it follows that $\max(P_{\Gamma}^{'})=\lbrace e_{r-1},e_r\rbrace$.

\subsection{Orbits of the adjoint representation and the isotropy algebra}

In the subsections before we have seen that the antipodal set in the tangent space consists of a union of orbits of the form $\Ad(K)\tx$ for $\tx\in \pi \max (\triangle^{'})$ or $\tx\in \pi \max (P_{\Gamma}^{'})$. It is well known that $\Ad(K)\textbf{x}$ is an embedded submanifold of $\p$, as it is the orbit of a Lie group action of a compact Lie group. Let $K_{\tx}$ be the stabilizer of $\tx$ under the adjoint action. The embedding $i$ has the form $i:K\slash K_{\tx}\to \Ad(K) \tx$. Let $\kk_{\tx}$ be the Lie algebra of $K_{\tx}$, which is called \emph{isotropy algebra}. Then\vskip -.6cm
\begin{align}\label{tangentspaceimmersion}
T_{\tx}\AT= di\: T_{eK_{\tx}} K\slash K_{\tx} = di\:\kk\slash\kk_{\tx}.
\end{align}
Hence the following proposition helps us to determine $T_{\tx}\AT$.

\bp (see \cite{Kon} proposition 2.1)\label{proposition-isotropy-algebra}\\
For a given $\tx\in \hp$, let $\Sigma_{\tx}:=\lbrace \alpha\in \Sigma\:|\:\alpha(\tx)=0\rbrace$. Then $\Sigma_{\tx}$ is a root subsystem. The isotropy subalgebra of the adjoint representation at $\tx$ is of the form $\kk_{\tx}=\kk(\Sigma_{\tx})$.
\ep

\section{The tangent space of the antipodal set}\label{section-tangent-space-of-antipodalset}

In this section we give an explicit description of the tangent space of the antipodal set of all irreducible symmetric spaces of compact type, but those with restricted root system $\mathfrak{a}_r$ and a non-trivial fundamental group different from $\mathbb{Z}_2$ or $\Z_{r+1}$. This enables us to calculate its dimensions in the next section.\\

We have seen in the previous section that, if we are not in one of the excluded cases, the antipodal set in the tangent space consists of disjoint orbits of $\tx$ with $\tx=\pi e_j$, $\tx=\frac{\pi}{2}(e_j+e_{j+1})$ for some specific $j$ or $\tx=\frac{\pi}{r+1}(e_1+\ldots +e_r)$. We fix one such $\tx$. Let $\mathbf{p}:=\exp_{eK}(\tx)$. By \eqref{antipodal-tangent-regular} and \eqref{tangentspaceimmersion} we have
\begin{align}\label{main-equation}
T_{\mathbf{p}}\AU=(d\exp_{eK})_{\tx} (T_{\tx}\AT)= (d\exp_{eK})_{\tx} di (\kk\slash \kk_{\tx}).
\end{align}

To calculate this tangent space we first need to determine $\kk\slash \kk_{\tx}$, then map it by $di$ to $\p$ and at last apply $(d\exp_{eK})_{\tx}$. 

For the first step Proposition \ref{proposition-isotropy-algebra}, \eqref{decompositionKandP} and \eqref{decomposition-subsystem} yield\vskip -.4cm

\begin{align*}
\kk\slash \kk_{\tx}=&\: \kk(0)\oplus\bigoplus_{\alpha\in \Sigma^{+}}\kk(\alpha)\slash(\kk(0)\oplus\bigoplus_{\alpha\in \Sigma_{\tx}^{+}}\kk(\alpha))\\
\cong & \bigoplus_{\alpha\in \Sigma^{+}\backslash\Sigma^{+}_{\tx}}\kk(\alpha).
\end{align*}

By definition $i$ maps $kK_{\tx}$ to $\Ad(k)\tx$, hence $(di)_{eK_\tx}=-\ad_{\tx}$. Let $Y\in \kk(\alpha)$, $\alpha\in\Sigma^{+}\backslash\Sigma^{+}_{\tx}$ and $H\in\hp$. It is easy to check that $\ad_H^2((di)_{eK_\textbf{x}}(Y)) = -\alpha(H)^2 [Y,\textbf{x}]$. This implies $(di)_{eK_{\textbf{x}}}(Y)\in\g(\alpha)$ and as $[Y,\textbf{x}]\in\p$, even $(di)_{eK_{\textbf{x}}}(Y)\in\p(\alpha)$. Since $di$ is of full rank and $\dim \kk(\alpha)=\dim\p(\alpha)$ (\cite{Lo} p. 60), this implies
\begin{align*}
(di)_{eK_{\textbf{x}}}(\kk(\alpha))=\p(\alpha).
\end{align*}

Let $p_0^{\tx}$ denote the parallel transport from $0$ to $\tx$. We summarize what we have just shown.
\bp\mbox{}\label{diisomorphism}\\
Notations as before, then
\begin{align*}
T_{\tx}\Ad(K)\tx=p_0^{\tx}(\bigoplus_{\alpha\in \Sigma^{+}\backslash\Sigma^{+}_{\tx}}\p(\alpha)).
\end{align*}
\ep

We determine now $\Sigma^{+}\backslash\Sigma^{+}_{\tx}$ for the three types of $\tx$ using Proposition \ref{proposition-isotropy-algebra}. Let $\lbrace\alpha_1,\ldots,\alpha_r\rbrace$ be the set of simple roots of $\Sigma$ with respect to the fixed ordering.\vsk
Let $\tx=\pi e_j$. Then $\Sigma_{\tx}=\lbrace \alpha=\sum_{k=1}^r c_k\alpha_k\:|\:\alpha(\pi e_j)=0\rbrace$. As $\alpha(\pi e_j)=0\Leftrightarrow c_j=0$, we get $\Sigma^{+}\backslash\Sigma^{+}_{\tx}= \lbrace \alpha=\sum_{k=1}^r c_k\alpha_k\:|\:c_j> 0\rbrace$.\vsk
Let $\tx=\frac{\pi}{2}(e_j+e_{j+1})$. With a similar argument we derive $\Sigma^{+}\backslash\Sigma^{+}_{\tx}= \lbrace \alpha=\sum_{k=1}^r c_k\alpha_k\:|\:c_j> 0 \vee c_{j+1}>0 \rbrace$.\vsk
Let $\tx=\frac{\pi}{r+1}(e_1+\ldots +e_r)$. Since $\alpha (\frac{\pi}{r+1}(e_1+\ldots +e_r))=\frac{\pi}{r+1}(c_1+\ldots +c_r)$ for $\alpha=\sum_{k=1}^r c_k\alpha_k$, we get $\Sigma_{\tx}=\emptyset$ and hence $\Sigma^{+}\backslash\Sigma^{+}_{\tx}=\Sigma^{+}$.\\

For the last step, namely applying $(d\exp_{eK})_{\tx}$, we use the following theorem.

\bt\label{kernelexpcor}(see \cite{Cr} p. 325 and Appendix)\\
Notations as before. Let $H\in\hp$, then the Euclidean parallel translate of 
\begin{align*}
\bigoplus\limits_{\alpha(H)\equiv 0\mod \pi\: ;\:\alpha(H)\neq 0\: ;\: \alpha\in\Sigma^{+}}\p (\alpha)\hspace{.2cm}
\end{align*}  
to $H$ constitutes the kernel of $d(\exp_{eK})_H$.
\et

We consider the three types of $\tx$ separately. Again $\psi=\sum_{i=1}^r d_i \alpha_i$ is the highest root of $\Sigma$ and $\alpha =\sum_{k=1}^{r} c_k \alpha_k \in \Sigma^{+}$ thus $c_k\in \mathbb{N}$.\vsk
Let $\tx=\pi e_j$. Then $\alpha(\tx)=\pi\alpha(e_j)= \pi\frac{c_j}{d_j}$. If we apply Theorem~\ref{kernelexpcor}, we see that $\p(\alpha)$ is part of the kernel of $(d\exp_{eK})_{\pi e_j}$ if and only if $\frac{c_j}{d_j}\in \mathbb{N}\backslash \lbrace 0\rbrace$. Thus
\begin{align*}
\ker((d\exp_{eK})_{\tx})=p_0^{\tx}\: (\:\bigoplus\limits_{\alpha \in J_j} \p(\alpha)\:),
\end{align*}
where $J_j:=\lbrace \alpha =\sum_{k=1}^{r} c_k \alpha_k \in \Sigma\:|\: \frac{c_j}{d_j}\in \mathbb{N} \backslash \lbrace 0\rbrace\rbrace$.\vsk
Let $\tx=\frac{\pi}{2}(e_j+e_{j+1})$. Then $\alpha(\tx)=\frac{\pi}{2}\alpha(e_j+e_{j+1})= \frac{\pi}{2}(\frac{c_j}{d_j}+\frac{c_{j+1}}{d_{j+1}})$. We see that $\p(\alpha)$ is part of the kernel of $(d\exp_{eK})_{\pi e_j}$ if and only if $\frac{c_j}{d_j}+\frac{c_{j+1}}{d_{j+1}}\in 2\mathbb{N} \backslash \lbrace 0\rbrace$. Thus
\begin{align*}
\ker((d\exp_{eK})_{\tx})=p_0^{\tx}\: (\:\bigoplus\limits_{\alpha \in J_{j,j+1}} \p(\alpha)\:),
\end{align*}
where $J_{j,j+1}:=\lbrace \alpha =\sum_{k=1}^{r} c_k \alpha_k \in \Sigma\:|\: \frac{c_j}{d_j}+\frac{c_{j+1}}{d_{j+1}}\in 2\mathbb{N} \backslash \lbrace 0\rbrace\rbrace$.\vsk

Let $\tx=\frac{\pi}{r+1}(e_1+\ldots +e_r)$. Then $\alpha (\frac{\pi}{r+1}(e_1+\ldots e_r))=\frac{\pi}{r+1}(c_1+\ldots +c_r)$. Note that we can assume in this case $\Sigma=\mathfrak{a}_r$ and therefore $c_i=0,1$ if $\alpha$ is a positive root. It follows $\frac{1}{r+1}(c_1+\ldots +c_r)\notin \mathbb{N}$. Hence 

\begin{align*}
\ker((d\exp_{eK})_{\tx})=\lbrace 0\rbrace.
\end{align*}

Let for the moment be $p:=p_o^{\tx}$. In view of \eqref{main-equation}, Proposition \ref{diisomorphism} and the results on $\ker((d\exp_{eK})_{\tx})$ we have just shown, we get 

\begin{align*}
T_{\exp_{eK}(\tx)}\exp_{eK}(\Ad(K)\textbf{x}) = & p\textbf{(}\bigoplus_{\alpha\in \Sigma^{+}\backslash\Sigma^{+}_{\tx}} \p(\alpha) \textbf{)}\slash p\textbf{(} \bigoplus\limits_{\alpha \in I} \p(\alpha) \cap \bigoplus_{\alpha\in \Sigma^{+}\backslash\Sigma^{+}_{\tx}}\p(\alpha)\textbf{)}\\
= & p (\bigoplus\limits_{\alpha \in I^{'}} \p(\alpha)),
\end{align*}
where $I=J_j, J_{j,j+1}$ or $\emptyset$, depending on which $\tx$ we are considering, and $I^{'}:=(\Sigma^{+}\backslash\Sigma^{+}_{\tx})\backslash (I \cap \Sigma^{+}\backslash\Sigma^{+}_{\tx})$. We summarize our results in a theorem.

\bt\label{non sc tangent space theorem}\mbox{}\\
Let $\G$ be an irreducible symmetric space of compact type, that has not both the restricted root system $\Sigma=\mathfrak{a}_r$ and a non-trivial $\Gamma\neq \Z_2$ or $\Z_{r+1}$, where $\Gamma$ is given through a Clifford-Klein-form of $\G$. Let $\psi=\sum_{i=1}^r d_i \alpha_i$ be the highest root of the restricted root system and assume roots are written in the way $\alpha =\sum_{k=1}^{r} c_k \alpha_k$. Furthermore let $\lbrace e_1,\ldots,e_r\rbrace$ be the corners of the Cartan polyhedron. If $\G$ is simply connected, then take $\tx\in \pi \max (\triangle^{'})$, if $\G$ is non-simply connected, then take $\tx\in \pi \max (P_{\Gamma}^{'})$, while $\max (\triangle^{'})$ and $\max (P_{\Gamma}^{'})$ are given in the tables \ref{table-maximal-corners-and-factors}, \ref{table-non-cs-max corners}. Then $\exp_{eK}(\tx)\in \AU$ and the tangent space of the antipodal set at $\exp_{eK}(\tx)$ can be given explicitly depending on the form of $\tx$. If $\tx=\pi e_j$, then 
\begin{align*}
T_{\exp_{eK}(\tx)} \AU = p_0^{\tx}( \bigoplus\limits_{\alpha \in J_j^{'}} \p(\alpha)), \spa \text{for } J_j^{'}=\lbrace \alpha \in \Sigma^{+}\:|\: \frac{c_j}{d_j}\notin \mathbb{N} \rbrace.
\end{align*}
If $\tx=\frac{\pi}{2}(e_j+e_{j+1})$, then
\begin{align*}
T_{\exp_{eK}(\tx)} \AU = p_0^{\tx}( \bigoplus\limits_{\alpha \in J_{j,j+1}^{'}} \p(\alpha)), \spa J_{j,j+1}^{'}=\lbrace \alpha \in \Sigma^{+}\:|\: \frac{c_j}{d_j}+\frac{c_{j+1}}{d_{j+1}}\notin 2\mathbb{N} \rbrace.
\end{align*}
If $\tx=\frac{\pi}{r+1}(e_1+\ldots +e_r)$, then 
\begin{align*}
T_{\exp_{eK}(\tx)} \AU= p_0^{\tx}( \bigoplus\limits_{\alpha \in \Sigma^{+}} \p(\alpha)).
\end{align*}
\et 

\begin{proof}
Almost everything has been proved so far. What is left, is to verify that $J_j^{'}=(\Sigma^{+}\backslash\Sigma^{+}_{\tx})\backslash (J_j \cap \Sigma^{+}\backslash\Sigma^{+}_{\tx}) $ and $J_{j,j+1}^{'}=(\Sigma^{+}\backslash\Sigma^{+}_{\tx})\backslash (J_{j,j+1} \cap \Sigma^{+}\backslash\Sigma^{+}_{\tx})$, which is straight forward.
\end{proof}

We know that for a positive root $\alpha =\sum_{k=1}^{r} c_k \alpha_k$ it is $c_k\in\mathbb{N}$. Let again $\tx=\pi e_j$. If the corresponding factor of the highest root $d_j$ equals 1, we have $c_j\slash d_j=c_j\in\mathbb{N}$. In this case $J_j^{'}=\emptyset$ and hence $\dim(T_{\exp_{eK}(\tx)}\AU)=0$. Thus some cases with $\dim(T_{\exp_{eK}(\tx)}\AU)=0$ can be read of tables \ref{table-maximal-corners-and-factors} and \ref{table-non-cs-max corners}. The calculations in the next section show that this are all cases. For simply connected $\G$ this result was first proved by Deng and Liu in \cite{DL}. We state it in the corollary below.

\bc\label{Deng-Liu-cor} (see \cite{DL})\mbox{}\\
Let $G\slash K$ be an irreducible compact simply connected symmetric space with restricted root system $\Sigma$. Let $\psi=\sum_{i=1}^r d_i\alpha_i$ be the highest root and let $e_j$ be a maximal corner of the Cartan polyhedron. If $d_j=1$, then $\dim \exp_{eK} (\Ad(K)(\pi e_j))=0$. In particular $\dim(\AU)=0$ if $\Sigma$ is one of the following: $\mathfrak{b}_2, \mathfrak{b}_3, \mathfrak{a}_r,\mathfrak{c}_r, \mathfrak{d}_r,\mathfrak{e_6}$ or $\mathfrak{e}_7$.
\ec

\section{Dimensions of the orbits building the antipodal set}\label{section-dimensions}

Let $\tx\in \pi\max(\triangle^{'})$ or $\tx\in \pi\max(P_{\Gamma}^{'})$ and $\tx\neq\frac{\pi}{r+1}(e_1+\ldots+e_r)$. In view of Theorem \ref{non sc tangent space theorem}, we have seen that $\tx=\pi e_j$ or $\tx=\frac{\pi}{2} (e_j+e_{j+1})$ and
\begin{align*}
\dim T_{\exp_{eK}(\tx)}\AU &= \dim \bigoplus\limits_{\alpha \in J_j^{'}} \p(\alpha),\hspace{.2cm} J_j^{'}=\lbrace \alpha =\sum_{k=1}^{r} c_k \alpha_k \in \Sigma^{+}\:|\: \frac{c_j}{d_j}\notin \mathbb{N} \rbrace\text{ or}\\
\dim T_{\exp_{eK}(\tx)} \AU &= \dim \bigoplus\limits_{\alpha \in J_{j,j+1}^{'}} \p(\alpha), \hspace{.2cm} J_{j,j+1}^{'}=\lbrace \alpha \:|\: \frac{c_j}{d_j}+\frac{c_{j+1}}{d_{j+1}}\notin 2\mathbb{N} \rbrace
\end{align*}
To determine the dimension of $\AU$ at $\exp_{eK}(\tx)$ we do the following: 
\begin{enumerate}
\item Determine the set $J^{'}_j$, if $\tx=\pi e_j$ or $J_{j,j+1}^{'}$, if $\tx=\frac{\pi}{2} (e_j+e_{j+1})$.
\item Sum up the dimensions of the spaces $\p(\alpha)$ for $\alpha\in J^{'}_j$ or $J_{j,j+1}^{'}$.
\end{enumerate}
If there are more than one maximal corners, we consider each corner and its orbit separately. The union of these orbits builds the antipodal set. The maximal corners and the corresponding factors $d_j$ are given in table~\ref{table-maximal-corners-and-factors} and table \ref{table-non-cs-max corners}. The facts on root systems can for example be found in \cite{Bou2}. The root multiplicities are well known, see for example \cite{Lo} p. 119 and p. 146. Throughout the rest of the section let $\lbrace x_1,\ldots,x_r \rbrace$ be a possibly scaled standard basis.\vskip .4cm
We give several calculations as example, as the other calculation are in the same manner we leave them away and just state the whole results in the tables at the end of this section.\\ 

\textbf{Simply connected symmetric spaces:}\\

Let $\G$ be a symmetric space of compact type with restricted root system $\Sigma=\mathfrak{bc}_r$. Then $\G$ is simply connected. A choice of simple roots of $\mathfrak{bc}_r$ is $\lbrace\alpha_1,\ldots,\alpha_{r-1},\alpha_r\rbrace=\lbrace x_1 -x_2,\ldots,x_{r-1}-x_r,x_r\rbrace$. Thus the positive roots are $\Sigma^{+}=\lbrace x_i \pm x_j (i<j), x_i, 2x_i\rbrace$. The maximal corner of the Cartan polyhedron is $e_r$, hence the antipodal set consists of one orbit and is therefore a manifold. The highest root is given by $\psi=2\sum_{j=1}^r \alpha_j$, meaning $d_r=2$. 
Thus $\frac{c_r}{d_r}\notin \mathbb{N}$ if and only if $c_r=1$. It is straight forward to check $J_r^{'}=\lbrace x_1,\ldots,x_r\rbrace$.\vskip .5cm

\textbf{Type A III with $\Sigma=\mathfrak{bc}_r$:} Given a symmetric space of type A III of the form $\G=SU(2r+q)\slash S(U(r) \times U(r+q))$, then the multiplicities of the roots $\alpha= \pm x_i$ are $m_{\alpha}=2q$. Hence $\dim \AU= 2qr$.\\

\textbf{Type C II with $\Sigma=\mathfrak{bc}_r$:} For a symmetric space of type C II of the form $G\slash K=SP(2r+q)\slash (SP(r)\times SP(r+q))$ the multiplicities of the roots $\alpha= \pm x_i$ are $m_{\alpha}=4q$. Thus $\dim \AU= 4qr$.\\

\textbf{Type D III with $\Sigma=\mathfrak{bc}_r$:} Let us consider a symmetric space of type D III. In this case the multiplicities of the roots $\alpha= \pm x_i$ are $m_{\alpha}=4$. Hence the dimension of the antipodal set of $G\slash K=SO(4r+2)\slash U(2r+1)$ takes the value $\dim \AU= 4r$.\\

\textbf{Type E III with $\Sigma=\mathfrak{bc}_2$:} For a symmetric space of type E III the multiplicities of the roots $\alpha= \pm x_i$ are $m_{\alpha}=8$. Hence, if $G\slash K=E_6\slash Spin(10)\cdot SO(2)$, then $\dim \AU= 16$.\\

\textbf{Type F II with $\Sigma=\mathfrak{bc}_1$:} A symmetric space of type F II  has root multiplicities $m_{\alpha}=8$ for $\alpha= \pm x_i$ . Thus for $G\slash K=F_4\slash Spin(9)$ we have $\dim \AU= 8$.\\

Let now $\G$ be a compact simply connected symmetric space with $\Sigma=\mathfrak{b}_r$. By table \ref{table-maximal-corners-and-factors} we get three different cases. For $r=2,3$ the antipodal set is of dimension zero (see Corollary \ref{Deng-Liu-cor}), for $r=4$ it is the union of a positive dimensional manifold and a point and for $r>4$ it is a positive dimensional manifold. In the following we determine the positive dimensional part, thus $r\geq 4$. A choice of simple roots is $\lbrace\alpha_1,\ldots,\alpha_{r-1},\alpha_r\rbrace=\lbrace x_1 -x_2,\ldots,x_{r-1}-x_r,x_r\rbrace$, then $\Sigma^{+}=\lbrace x_i \pm x_j\: (i<j), x_i\rbrace$. The maximal corner of the Cartan polyhedron is $e_r$ with $d_r=2$. We get $J_r^{'}=\lbrace x_1,\ldots, x_r \rbrace$.\\

\textbf{Type BD I with $\Sigma=\mathfrak{b}_r$:} For a simply connected symmetric space of type BD I the multiplicities of the roots $\alpha= \pm x_i$ are $m_{\alpha}=q$. Thus for $G\slash K=G_{r,r+q}$ with $r\geq 4$ and $q>0$, the positive dimensional part of the antipodal set is of dimension $rq$.\\

\textbf{Spin(2r+1):} We consider the Lie group $Spin(2r+1)$ with $r\geq 4$. In this case the multiplicities of all roots are $m_{\alpha}=2$ and thus the positive dimensional part of the antipodal set is of dimension $2r$.\\

\textbf{Non-simply connected symmetric spaces:}\\

\textbf{Type BD I with $\Sigma=\mathfrak{b}_r$ and $\Gamma=\mathbb{Z}_2$:}
Let $\tM\slash\Gamma$ be an irreducible symmetric space given in the Clifford-Klein-form such that $\tM=Gr_{r,r+q}$, $r\geq 2$, $q>0$, $\Sigma=\mathfrak{b}_r$ and $\Gamma=\mathbb{Z}_2$. In this case $\max(P_{\Gamma}^{'})=\lbrace e_r\rbrace$ and $d_r=2$. A choice of simple roots is $\lbrace\alpha_1,\ldots,\alpha_{r-1},\alpha_r\rbrace=\lbrace x_1 -x_2,\ldots,x_{r-1}-x_r,x_r\rbrace$, hence $\Sigma^{+}=\lbrace x_i \pm x_j\: (i<j), x_i\rbrace$. It is easy to check that $J_r^{'}=\lbrace x_1,\ldots, x_r \rbrace$. The root multiplicities of those roots are $q$, which implies $\dim\mathrm{A}(eK)= rq$.\\

\textbf{Type A I with $\Gamma=\mathbb{Z}_2$:}
Let $\tilde{M}\slash\Gamma$ be an irreducible symmetric space of type A I such that $\Gamma=\mathbb{Z}_2$. In this case $\Sigma=\mathfrak{a}_r$.\\ 
If $r\geq 3$, $r$ odd and $\frac{r+1}{2}$ even, then $\max(P_{\Gamma}^{'})=\lbrace e_{\frac{r+1}{4}} \rbrace$. As $d_{\frac{r+1}{4}}=1$, it follows $\dim \mathrm{A}(eK)=0$.\\ 
If $r\geq 3$, $r$ and $\frac{r+1}{2}$ odd, then $\max(P_{\Gamma}^{'})=\lbrace \frac{1}{2} (e_{\frac{r-1}{4}}+e_{\frac{r+3}{4}})\rbrace$ and $d_{\frac{r-1}{4}}=1=d_{\frac{r+3}{4}}$. Since $\Sigma^{+}=\lbrace x_i - x_j\:|\: i<j\leq r+1\rbrace$ and we are looking for those roots with $c_{\frac{r-1}{4}}+c_{\frac{r+3}{4}}=1$, we get $J_{\frac{r-1}{4},\frac{r+3}{4}}^{'}=\lbrace x_i - x_{\frac{r+3}{4}} \:|\: i < \frac{r+3}{4}\rbrace \cup \lbrace x_{\frac{r+3}{4}} - x_j \:|\:  \frac{r+3}{4}<j\leq r+1 \rbrace$. As the root multiplicities of all roots are 1, we deduce $\dim \mathrm{A}(eK)=r$.\\

\textbf{Type C II with $\Gamma=\mathbb{Z}_2$:}
Let $\tilde{M}\slash\Gamma$ be an irreducible symmetric space of type C II such that $\Gamma=\mathbb{Z}_2$.  In this case $\Sigma=\mathfrak{c}_r$. We take the simple roots to be $\lbrace\alpha_1,\ldots,\alpha_{r-1},\alpha_r\rbrace=\lbrace x_1 -x_2,\ldots,x_{r-1}-x_r,2x_r\rbrace$, which yields $\Sigma^{+}=\lbrace x_i \pm x_j \:(i<j), 2x_i\rbrace$.\\ 
If $r$ is even, we have $\max(P_{\Gamma}^{'})=\lbrace e_{\frac{r}{2}} \rbrace$ and $d_{\frac{r}{2}}=2$. We are looking for those roots that have $c_{\frac{r}{2}}=1$. The result is $J_{\frac{r}{2}}^{'}=\lbrace x_i \pm x_j\: |\: i\leq \frac{r}{2}<j \rbrace$. As the root multiplicities of $x_i \pm x_j$ with $i\neq j$ are 4, we derive $\dim \mathrm{A}(eK)=2r^2$.\\
If $r$ is odd, we have $\max(P_{\Gamma}^{'})=\lbrace \frac{1}{2} (e_{\frac{r-1}{2}}+e_{\frac{r+1}{2}}) \rbrace$ and $d_{\frac{r-1}{2}}=2=d_{\frac{r+1}{2}}$. Then $\alpha=\sum_{i=1}^r c_i\alpha_i\in J_{\frac{r-1}{2},\frac{r+1}{2}}^{'}$ if and only if $c_{\frac{r-1}{2}} + c_{\frac{r+1}{2}}= 1,2$ or $3$. Thus 
\[J_{\frac{r-1}{2},\frac{r+1}{2}}^{'}= \lbrace x_i \pm x_j \:|\: i\leq \frac{r+1}{2}, i<j,j\geq \frac{r+1}{2} \rbrace\cup \lbrace 2x_{\frac{r+1}{2}}\rbrace.\] 
The root multiplicities are $m_{x_i \pm x_j}=4$ for $i\neq j$ and $m_{2x_i}=3$ and hence $\dim \mathrm{A}(eK)=2r^2+4r -3$.\\

\textbf{Spin(2r)$\slash \Gamma$ with $|\Gamma|=4$:}
Consider the Lie group $Spin(2r)\slash \Gamma$ with $|\Gamma|=4$. Then $\Sigma=\mathfrak{d}_r$. A set of simple roots can be chosen to be $\lbrace x_1 - x_2,\ldots, x_{r-1}-x_r, x_{r-1}+ x_r\rbrace$. In this case we get $\Sigma^{+}=\lbrace x_i\pm x_j\:|\: i<j\rbrace$.\\
If $r$ is even, then $\Gamma=\mathbb{Z}_2\oplus\mathbb{Z}_2$. The maximal corner is $e_{\frac{r}{2}}$ and $d_{\frac{r}{2}}=2$. One can check that $J_{\frac{r}{2}}^{'}= \lbrace x_i \pm x_j\:|\: i\leq \frac{r}{2}< j\rbrace$. The root multiplicities are 2, thus $\dim\mathrm{A}(eK)=r^2$.\\
If $r$ is odd, then $\Gamma=\mathbb{Z}_4$. The maximal corner is $\frac{1}{2} (e_{\frac{r-1}{2}}+e_{\frac{r+1}{2}})$ and $d_{\frac{r-1}{2}}=2=d_{\frac{r+1}{2}}$. Similar as before we are looking for roots with $c_{\frac{r-1}{2}}+c_{\frac{r+1}{2}}=1,2$ or $3$. We get $J_{\frac{r-1}{2},\frac{r+1}{2}}^{'}= \lbrace x_i \pm x_j \:|\: i<j, i\leq \frac{r+1}{2},j\geq \frac{r+1}{2} \rbrace$. The root multiplicities are 2, thus $\dim\mathrm{A}(eK)=r^2+2r-3$.\\

\textbf{Type E I with $\Gamma=\Z_3$:} A choice of simple roots is 
\begin{align*}
\alpha_1&=\frac{1}{2}(x_1+x_8-x_2-x_3\ldots -x_7), \quad &&\alpha_2=x_1+x_2, \quad &&&\alpha_3=x_2-x_1,\\
\alpha_4 & = x_3 - x_2, \quad &&\alpha_5=x_4-x_3 \quad &&&\alpha_6=x_5 - x_4.
\end{align*} 
In this case
\begin{align*}
\Sigma^{+}= \hspace{.5cm}& \lbrace x_i\pm x_j \:|\: 1\leq j<i\leq 5\rbrace\\
\cup \hspace{.1cm}& \lbrace \frac{1}{2}(x_8-x_7-x_6+\sum\limits_{i=1}^{5} (-1)^{v(i)}x_i)\:|\: \sum\limits_{i=1}^{5}v(i) \text{ is even}\rbrace.
\end{align*}

The maximal corner is $e_4$ and $d_4=3$. Thus we are looking for roots such that $c_4=1$ or $2$. There are no roots with $c_4>3$. We get the following result
\begin{align*}
J_4^{'}= & \hspace{.4cm}\lbrace x_i - x_j\:|\: j\leq 2<i\leq 5\rbrace \hspace{.68cm}\cup \lbrace x_i +x_j\:|\:j<i, 3\leq i\leq 5 \rbrace\\
& \cup \lbrace \alpha_1+x_i -x_1 \:|\: i=3,4 \text{ or } 5\rbrace \cup \lbrace \alpha_1 +x_i+x_j\:|\:2\leq j<i\leq 5\rbrace\\
& \cup \lbrace \alpha_1+\alpha_3+x_i+x_j\:|\: (j,i)=(3,4),(3,5) \text{ or } (4,5)\rbrace
\end{align*}

As $|J_4^{'}|=27$ and $m_{\alpha}=1$ for all roots, we get $\dim\mathrm{A}(eK)=27$.
\vskip .8cm

Let us consider some cases, where $\tx=\frac{\pi}{r+1}(e_1+\ldots+e_r)\in \max(P^{'}_{\Z_{r+1}})$. If $\tx$ is of that form, then it has to be $\Sigma=\mathfrak{a}_r$. It is well known that in this case all root multiplicities are the same. As now $T_{\exp_{eK}(\tx)} \AU= p_0^{\tx}( \bigoplus_{\alpha \in \Sigma^{+}} \p(\alpha))$, we get $\dim\mathrm{A}(eK)=m_{\alpha}|\mathfrak{a}^{+}_r|$, where $m_{\alpha}$ is the root multiplicity and $|\mathfrak{a}^{+}_r|=\frac{r(r+1)}{2}$.\vsk
\textbf{Type A I with rank $r$ and $\Gamma=\Z_{r+1}$:} We have $m_{\alpha}=1$, therefore $\dim\mathrm{A}(eK)=\frac{r(r+1)}{2}$.\vsk
\textbf{Type A II with rank $r$ and $\Gamma=\Z_{r+1}$:} We have $m_{\alpha}=4$, therefore $\dim\mathrm{A}(eK)=2r(r+1)$.\vsk
\textbf{Type E IV with $\Gamma=\Z_{3}$:} We have $m_{\alpha}=8$, this gives $\dim\mathrm{A}(eK)=24$.\vsk
\textbf{SU(r+1) with $\Gamma=\Z_{r+1}$:} In this case $m_{\alpha}=2$, therefore $\dim\mathrm{A}(eK)=r(r+1)$.\\

The following four tables contain the dimensions of the different components of the antipodal set of all irreducible symmetric spaces of compact type, but those cases explicitly excluded. We want to give a few remarks on those tables: If there is more than one maximal corner, then the dimensions of the corresponding orbits in the following column are listed in the same order as the maximal corners. The indexing of the simple roots is again as in \cite{Hel} p. 477, 478. In the following $q$ is a positive integer. Furthermore we use the notation $e=eK$ and $p_i=\exp_{eK}(\pi e_i)$. Other references were given at the beginning of this section. 

\newpage

\begin{center}
\captionof{table}{Maximal corners of the Cartan polyhedron and the dimensions of the components of the antipodal set of irreducible compact simply connected symmetric spaces of type I} \label{table-type-I-dimension} 
\setlength{\extrarowheight}{2.2pt}
\begin{tabular}{|p{1.1cm}|p{6cm}|p{.9cm}|p{1.4cm}|p{1.1cm}|}
\hline
Type & $M$ or $(\g,\kk)$ & $\Sigma$ & $\max(\triangle^{'})$ & $\dim\mathrm{A}$\\  \hline
A I & $SU(2r)\slash SO(2r)$ & $\mathfrak{a}_{2r-1}$ & $e_{r}$ & 0  \\ 
& $SU(2r+1)\slash SO(2r+1)$  & $\mathfrak{a}_{2r}$ & $e_{r}$; $e_{r+1}$ & 0; 0  \\ \hline
A II & $SU(4r)\slash Sp(2r)$ & $\mathfrak{a}_{2r-1}$ & $e_{r}$ & 0  \\ 
& $SU(4r+2)\slash Sp(2r+1)$ & $\mathfrak{a}_{2r}$ & $e_{r}$; $e_{r+1}$ & 0; 0  \\ \hline
A III & $Gr_{r,r+q}(\mathbb{C})$,\spa$r\geq 2, q\geq 1$ or $r=1$  & $\mathfrak{bc}_{r}$ & $e_r$ & $2qr$\\ 
 & $Gr_{r,r}(\mathbb{C})$,\spa\spa\spa $r\geq 2$ & $\mathfrak{c}_{r}$ & $e_r$ & 0\\ \hline

C I & $Sp(r)\slash U(r)$ &  $\mathfrak{c}_{r}$ & $e_r$ & 0 \\  \hline

C II & $Gr_{r,r+q}(\mathbb{H})$,\spa$r\geq 2, q\geq 1$ or $r=1$  & $\mathfrak{bc}_{r}$ & $e_r$ & $4qr$ \\  
& $Gr_{r,2r}(\mathbb{H})$,\spa\spa $r\geq 2$ & $\mathfrak{c}_{r}$ & $e_r$ & 0\\ \hline
 
BD I & $Gr_{r,r+q}$, \spa $r=2,3$ & $\mathfrak{b}_{r}$ & $e_1$ & 0 \\
& $Gr_{4,4+q}$  & $\mathfrak{b}_{4}$ & $e_1$; $e_4$ & 0; $4q$\\
& $Gr_{r,r+q}$, \spa $r\geq 5$ & $\mathfrak{b}_{r}$ & $e_r$ & $rq$\\  
& $Gr_{1,1+q}$  & $\mathfrak{a}_1$ & $e_1$ & 0 \\  
& $Gr_{4,8}$ & $\mathfrak{d}_{4}$ & $e_1;e_3;e_4$ & 0; 0; 0 \\ 
& $Gr_{r,2r}$, \spa $r\geq 5$ & $\mathfrak{d}_{r}$ & $e_{r-1};e_r$ & 0; 0\\  \hline
 
D III & $SO(4r)\slash U(2r)$  & $\mathfrak{c}_{r}$ & $e_{r}$ & 0 \\
& $SO(4r+2)\slash U(2r+1)$  & $\mathfrak{bc}_{r}$ & $e_{r}$ & $4r$  \\ \hline
E I & $(\mathfrak{e}_6,\mathfrak{sp}(4))$ & $\mathfrak{e}_{6}$ & $e_1$; $e_6$ & 0; 0 \\ \hline
E II & $(\mathfrak{e}_6,\mathfrak{su}(6)\oplus\mathfrak{su}(2))$ & $\mathfrak{f}_{4}$ & $e_4$ & 16 \\ \hline
E III & $(\mathfrak{e}_6, \mathfrak{so}(10)\oplus\R)$  & $\mathfrak{bc}_{2}$ & $e_2$ & 16 \\ \hline
E IV & $(\mathfrak{e}_6,\mathfrak{f}_4)$  & $\mathfrak{a}_{2}$ & $e_1$; $e_2$ & 0; 0 \\ \hline
E V & $(\mathfrak{e}_7,\mathfrak{su}(8))$  & $\mathfrak{e}_{7}$ & $e_7$ & 0 \\ \hline
E VI & $(\mathfrak{e}_7,\mathfrak{so}(12)\oplus \mathfrak{su}(2))$ & $\mathfrak{f}_{4}$ & $e_4$ & 32 \\ \hline
E VII & $(\mathfrak{e}_7,\mathfrak{e}_6\oplus \R)$ &  $\mathfrak{c}_{3}$ & $e_3$ & 0 \\ \hline
E VIII & $(\mathfrak{e}_8,\mathfrak{so}(16))$ & $\mathfrak{e}_{8}$ & $e_1$ & 64 \\ \hline
E IX & $(\mathfrak{e}_8,\mathfrak{e}_7\oplus\mathfrak{su}(2))$ &  $\mathfrak{f}_{4}$ & $e_4$ & 64 \\ \hline
F I & $(\mathfrak{f}_4,\mathfrak{sp}(3)\oplus\mathfrak{su}(2)) $ &  $\mathfrak{f}_{4}$ & $e_4$ & 8 \\ \hline
F II & $(\mathfrak{f}_4,\mathfrak{so}(9))$ &  $\mathfrak{bc}_{1}$ & $e_1$ & 8 \\ \hline
G & $(\mathfrak{g}_2,\mathfrak{su}(2)\oplus\mathfrak{su}(2))$ & $\mathfrak{g}_{2}$ & $e_1$ & 3 \\ \hline
\end{tabular}\\
\end{center}

\newpage
\begin{center}
\captionof{table}{Dimensions of the components of the antipodal set of irreducible symmetric spaces of compact type and type I} \label{table-non-sc-type-I-dimension} 
\setlength{\extrarowheight}{2.2pt}
\begin{tabular}{|p{.95cm}|p{3.9cm}|p{0.8cm}|p{3.1cm}|p{2.1cm}|}
\hline
Type & $\tM$ or $(\g,\kk)$ & $\Sigma$ & $\Gamma$ & $\dim\mathrm{A}$\\  \hline
A I & $SU(2r+2)\slash SO(2r+2)$, $\frac{r+1}{2}$ even  & $\mathfrak{a}_{2r+1}$ & $\mathbb{Z}_2$ & 0   \\
& $SU(2r+2)\slash SO(2r+2)$, $\frac{r+1}{2}$ odd  & $\mathfrak{a}_{2r+1}$ & $\mathbb{Z}_2$ & $2r +1$  \\ 
& $SU(r+1)\slash SO(r+1)$ & $\mathfrak{a}_{r}$ & $\mathbb{Z}_{r+1}$ & $\frac{r(r+1)}{2}$  \\ 
& otherwise & & & unknown\\ \hline
A II & $SU(4r+4)\slash Sp(2r+2)$, $\frac{r+1}{2}$ even  & $\mathfrak{a}_{2r+1}$ & $\mathbb{Z}_2$ & 0 \\
& $SU(4r+4)\slash Sp(2r+2)$, $\frac{r+1}{2}$ odd  & $\mathfrak{a}_{2r+1}$ & $\mathbb{Z}_2$ & $8r+4$ \\ 
& $SU(2r+2)\slash Sp(r+1)$ & $\mathfrak{a}_{r}$ & $\mathbb{Z}_{r+1}$ & $2r(r+1)$ \\ 
& otherwise & &  & unknown\\ \hline

A III  & $Gr_{r,r}(\mathbb{C})$, \hspace{.1cm} $r\geq 2$, r even & $\mathfrak{c}_{r}$ & $\mathbb{Z}_2$ & $r^2$\\
& $Gr_{r,r}(\mathbb{C})$, \hspace{.1cm} $r\geq 2$, r odd & $\mathfrak{c}_{r}$ & $\mathbb{Z}_2$ & $r^2+2r-2$\\ \hline

C I & $Sp(r)\slash U(r)$ \hspace{.2cm} r even &  $\mathfrak{c}_{r}$ & $\Z_2$ & $\frac{1}{2}r^2$ \\ 
& $Sp(r)\slash U(r)$ \hspace{.2cm} r odd &  $\mathfrak{c}_{r}$ & $\Z_2$ & $\frac{1}{2}(r^2+2r-1)$ \\ \hline

C II & $Gr_{r,r}(\mathbb{H})$,\hspace{.1cm} $r\geq 2$, r even & $\mathfrak{c}_{r}$ & $\Z_2$ & $2r^2$\\ 
& $Gr_{r,r}(\mathbb{H})$,\hspace{.1cm} $r\geq 2$, r odd & $\mathfrak{c}_{r}$ & $\Z_2$ & $2r^2+4r-3$ \\ \hline
 
BD I & $Gr_{r,r+q}$\spa $r\geq 2$ & $\mathfrak{b}_{r}$ & $\Z_2$ & $rq$ \\  
& $Gr_{r,r}$ \spa\spa $r$ even & $\mathfrak{d}_{r}$ & $\Z_2\oplus\Z_2$ & $\frac{1}{2}r^2$ \\ 
& $Gr_{r,r}$ \spa\spa $r$ odd & $\mathfrak{d}_{r}$ & $\Z_4$ & $\frac{1}{2}(r^2+2r-3)$ \\ 
& $Gr_{r,r}$ \spa\spa $r$ even & $\mathfrak{d}_{r}$ & $\lbrace e,p_1\rbrace$ & 0; 0\\ 
& $Gr_{r,r}$ \spa\spa $r\leq 6$, $r$ even & $\mathfrak{d}_{r}$ & $\lbrace e,p_{r-1}\rbrace$ or $\lbrace e,p_r\rbrace$  & 0\\ 
& $Gr_{8,8}$  & $\mathfrak{d}_{8}$ & $\lbrace e,p_{r-1}\rbrace$ or $\lbrace e,p_r\rbrace$ & 0; $\frac{1}{2}r^2$ \\ 
& $Gr_{r,r}$ \spa\spa $r\geq 10$, $r$ even & $\mathfrak{d}_{r}$ & $\lbrace e,p_{r-1}\rbrace$ or $\lbrace e,p_r\rbrace$ & $\frac{1}{2}r^2$ \\  \hline

D III & $SO(4r)\slash U(2r)$, \spa $r$ even & $\mathfrak{c}_{r}$ & $\Z_2$ & $2r^2$ \\ 
& $SO(4r)\slash U(2r)$, \spa $r$ odd & $\mathfrak{c}_{r}$ & $\Z_2$ & $2r^2+4r-5$ \\ \hline
E I & $(\mathfrak{e}_6,\mathfrak{sp}(4))$ & $\mathfrak{e}_{6}$ & $\Z_3$ & $27$ \\ \hline
E IV & $(\mathfrak{e}_6,\mathfrak{f}_4)$  & $\mathfrak{a}_{2}$ & $\Z_3$ & $24$ \\ \hline
E V & $(\mathfrak{e}_7,\mathfrak{su}(8))$  & $\mathfrak{e}_{7}$ & $\Z_2$ & $35$ \\ \hline
E VII & $(\mathfrak{e}_7,\mathfrak{e}_6\oplus\R)$ &  $\mathfrak{c}_{3}$ & $\Z_2$ & $49$ \\ \hline
\end{tabular}
\end{center}

\newpage
\captionof{table}{Maximal corners of the Cartan polyhedron and the dimensions of the components of the antipodal set of irreducible compact simply connected symmetric spaces of type II} \label{table-type-II-dimension} 
\setlength{\extrarowheight}{2.2pt}
\begin{tabular}{|p{5cm}|p{1cm}|p{2.3cm}|p{2.5cm}|} 
\hline
$G$ & $\Sigma$ & $\max(\triangle^{'})$ & $\dim\mathrm{A}$\\ \hline  
$SU(2r)$ & $\mathfrak{a}_{2r-1}$ & $e_{r}$ & 0 \\
$SU(2r+1)$ & $\mathfrak{a}_{2r}$ & $e_{r}$; $e_{r+1}$ & 0; 0 \\ \hline
$Spin(2r+1)$ \hspace{.2cm} $r=2,3$  & $\mathfrak{b}_{r}$ & $e_1$ & 0 \\ 
$Spin(9)$  & $\mathfrak{b}_{4}$ &$e_1$; $e_4$ & 0; 8\\ 
$Spin(2r+1)$ \hspace{.2cm} $r>4$  & $\mathfrak{b}_{r}$ & $e_r$ & 2r \\ \hline
$Sp(r)$ & $\mathfrak{c}_{r}$ & $e_r$ & 0 \\ \hline
$Spin(8)$ & $\mathfrak{d}_{4}$ & $e_1;e_3$; $e_4$ & 0; 0; 0 \\ 
$Spin(2r)$\hspace{.2cm} $r\geq 5$ & $\mathfrak{d}_{r}$ & $e_{r-1}$; $e_r$ & 0; 0\\ \hline
$E_6$ & $\mathfrak{e}_{6}$ & $e_1$; $e_6$ & 0; 0 \\ \hline
$E_7$ &  $\mathfrak{e}_7$ & $e_7$ & 0 \\ \hline
$E_8$ &  $\mathfrak{e}_{8}$ & $e_1$ & 128 \\ \hline
$F_4$ &  $\mathfrak{f}_{4}$ & $e_4$ & 16 \\ \hline
$G_2$ & $\mathfrak{g}_{2}$ & $e_1$ & 6 \\ \hline
\end{tabular}\vskip 1.2cm

\captionof{table}{Dimensions of the components of the antipodal set of irreducible symmetric spaces of compact type and type II} \label{table-non-cs-type-II-dimension} 
\setlength{\extrarowheight}{2.2pt}
\begin{tabular}{|p{4.2cm}|p{1cm}|p{3.1cm}|p{2.5cm}|}
\hline
$\tilde{G}$ & $\Sigma$ & $\Gamma$ & $\dim\mathrm{A}$\\ \hline  
$SU(2r+2)$ & $\mathfrak{a}_{2r+1}$ & $\Z_2$ & $4r+2$ \\
$SU(r+1)$ & $\mathfrak{a}_{r}$ & $\Z_{r+1}$ & $r(r+1)$ \\
otherwise &  &  & unknown \\ \hline
$Spin(2r+1)$  & $\mathfrak{b}_{r}$ & $\Z_2$ & $2r$ \\ \hline
$Sp(r)$ \spa $r$ even & $\mathfrak{c}_{r}$ & $\Z_2$ & $r^2$ \\
$Sp(r)$ \spa $r$ odd & $\mathfrak{c}_{r}$ & $\Z_2$ & $r^2+2r-1$\\ \hline
$Spin(2r)$ \spa $r$ even & $\mathfrak{d}_{r}$ & $\Z_2\oplus\Z_2$ & $r^2$ \\ 
$Spin(2r)$\hspace{.2cm} $r$ odd & $\mathfrak{d}_{r}$ & $\Z_4$ & $r^2+2r-3$\\
$Spin(2r)$\hspace{.2cm} $r$ even & $\mathfrak{d}_{r}$ & $\lbrace e, p_1\rbrace$ & 0; 0\\
$Spin(2r)$\hspace{.2cm} $r=4,6$ & $\mathfrak{d}_{r}$ & $\lbrace e,p_{r-1}\rbrace$ or $\lbrace e,p_r\rbrace$ & 0\\
$Spin(16)$ & $\mathfrak{d}_{8}$ & $\lbrace e,p_{r-1}\rbrace$ or $\lbrace e,p_r\rbrace$ & 0; $64$\\
$Spin(2r)$\hspace{.2cm} $r\geq 10$, $r$ even & $\mathfrak{d}_{r}$ & $\lbrace e,p_{r-1}\rbrace$ or $\lbrace e,p_r\rbrace$ & $r^2$\\ \hline
$E_6$ & $\mathfrak{e}_{6}$ & $\Z_3$ & 54 \\ \hline
$E_7$ &  $\mathfrak{e}_7$ & $\Z_2$ & 70 \\ \hline

\end{tabular}\vskip 1.2cm

\section{Appendix}

In the sections above we have cited a few results that are not stated in the literature exactly in the way we presented them. For Theorem \ref{Yang-theo}, Theorem \ref{Yang2-theorem} and Theorem \ref{kernelexpcor} this goes back to the fact that we use a different definition of the restricted root system. While the maximal corners we cited are not exactly given in the literature, they can be easily derived as we show below. 

\subsection*{Yangs theorems}

We want to derive Theorem \ref{Yang-theo} from Theorem 1.3. in \cite{Yang}. As Theorem \ref{Yang2-theorem} can be derived from Theorem 4.1. in \cite{Yang2} in a similar way, we leave that away.\\ 
Let $M=G\slash K$ be an irreducible compact simply connected symmetric space. Let $\g=Lie (G)$ and $\g^{\C}$ it's complexification. In view of \cite{Lo} p. 73 there is another way to define restricted roots, namely by the following:\\

\begin{large}
Definition:\\ \end{large} 
\textit{Let $\g^{\C}=\kk^{\C}+\p^{\C}$ be the decomposition coming from the Cartan decomposition $\g=\kk+\p$. Then $\h_{\p^{\C}}:=\hp^{\C}$ is a maximal abelian subspace in $\p^{\C}$. Furthermore let 
\begin{align*}
\g^{\C}_{\gamma}:=\lbrace X\in\g^{\C} \:|\: \ad_H (X)=\gamma(H)X \hspace{.4cm}\forall H\in\h_{\p^{\C}}\rbrace.
\end{align*}
Then $\gamma\in\hat{\Sigma}$ if and only if $\gamma\in\h_{\p^{\C}}^{*}$, such that $\g_{\gamma}\neq\lbrace 0\rbrace$ and $\gamma\neq 0$. The set $\hat{\Sigma}$ is again called restricted root system.}\\

This is an equivalent description to the one Yang uses in \cite{Yang} (see \cite{Lo} p. 73). It is well known that $\hat{\Sigma}$ is an irreducible abstract root system. Let $\Pi^{\hat{\Sigma}}$ be the simple roots of $\hat{\Sigma}$ and $\psi$ the highest root of $\hat{\Sigma}$. In the same manner as before in \eqref{Cartan-polyhedron} and \eqref{cartanpolyhedronedge}, let
\begin{align*}
\triangle_{\hat{\Sigma}}=\lbrace x\in \sqrt{-1}\h_{\p}\:|\: \gamma(x)\geq 0 \text{ for } \gamma \in \Pi^{\hat{\Sigma}} \wedge \psi(x)\leq 1 \rbrace,
\end{align*}
while $\hp$ and hence also $\sqrt{-1}\hp$ are naturally $\R$-subspaces of $\h_{\p^{\C}}$. The edges that do not contain 0 are denoted by\vskip -.4cm

\begin{align*}
\triangle^{'}_{\hat{\Sigma}}:=\lbrace x\in \triangle_{\hat{\Sigma}}\:|\: \psi(x)= 1\rbrace.
\end{align*}

Now we are able to cite the theorem by L. Yang.

\bt\label{CartanPMain} (see \cite{Yang} p. 689)\\
Let $M=G\slash K$ be a compact simply connected symmetric space. Then $C_{T}(eK)=\Ad(K)(\pi \sqrt{-1}\triangle^{'}_{\hat{\Sigma}})$.
\et

To get that version we wrote in Theorem \ref{Yang-theo} we need to show that $\triangle^{'}=\sqrt{-1}\triangle^{'}_{\hat{\Sigma}}$. Since $\triangle^{'}_{\hat{\Sigma}}\subset \sqrt{-1}\hp$ it is clear that $\sqrt{-1}\triangle^{'}_{\hat{\Sigma}}\subset \hp$. We use the following:

\bl (see \cite{Lo} p. 58, 60)\\
Let $G\slash K$ be a symmetric space and $\Sigma\subset\h_{\p}^{*}$, $\hat{\Sigma}\subset\h_{\p^{\C}}^{*}$ the corresponding restricted root systems described above. Then
\begin{align*}
\sqrt{-1}\hat{\Sigma}=\Sigma.
\end{align*}
\el

This lemma in particular implies that if $\gamma\in \hat{\Sigma}$, then $\sqrt{-1}\gamma\in \hp^{*}$. Having given a set of positive roots $\Sigma^{+}$, the simple roots are those positive roots that can not be written as the sum of two others. This implies that if we take a set of positive roots $\Sigma^{+}$ and put $\hat{\Sigma}^{+}:=\sqrt{-1}\Sigma^{+}$, then the sets of simple roots are related by $\Pi^{\Sigma}=\sqrt{-1}\Pi^{\hat{\Sigma}}$, while it is routine to show that $\sqrt{-1}\Sigma^{+}$ is a set of positive roots. We denote the Weyl chambers of $\Sigma^{+}$ and $\hat{\Sigma}^{+}$ by $C_{\Sigma^{+}}$ and $C_{\hat{\Sigma}^{+}}$, respectively. Let $\sqrt{-1} H\in C_{\hat{\Sigma}^{+}}\subset\sqrt{-1}\hp$ and $\gamma\in \Pi^{\hat{\Sigma}}$. Then
\begin{align*}
\gamma(\sqrt{-1}H)> 0 \quad \Longleftrightarrow \quad \sqrt{-1}\gamma (H)>0.
\end{align*}
As $\sqrt{-1}\gamma\in \Pi^{\Sigma}$, all elements in $\Pi^{\Sigma}$ can be won in this way, and $\gamma$ was an arbitrary simple root, it follows $H\in C_{\Sigma^{+}}$. With another argument of that form it follows $C_{\Sigma^{+}}=\sqrt{-1} C_{\hat{\Sigma}^{+}}$. Let $\psi\in \hat{\Sigma}$ such that $\sqrt{-1}\psi$ is the highest root of $\Sigma$, $\gamma\in\hat{\Sigma}$ and $H\in C_{\Sigma^{+}}$. Then
\begin{align*}
\sqrt{-1}\psi (H)\geq \sqrt{-1}\gamma(H) \quad \Longleftrightarrow\quad \psi (\sqrt{-1}H)\geq \gamma (\sqrt{-1}H).
\end{align*}
This implies that $\psi$ is the highest root of $\Sigma$. In particular it follows that $\triangle^{'}=\sqrt{-1}\triangle^{'}_{\hat{\Sigma}}$.

\subsection*{Maximal corners of irreducible non-simply connected symmetric spaces}

We shortly explain how the maximal corners of $P_{\Gamma}$ can be derived from \cite{Yang2}. For each of the 9 cases Yang has considered we give the approach to determine the maximal corners out of his results. 

In the cases II, VII, VIII and IX he determined all the corners of $P_{\Gamma}$, thus the results given in table \ref{table-non-cs-max corners} can be easily verified. Case VI was described above and case III is similar to that. In case I Yang proves that the there is only one maximal corner and it is straight forward to verify that the given one is the right one. Let's consider case IV and V. As the diameter of $P_{\Gamma}$ is given we can easily check that the given points in table \ref{table-non-cs-max corners} are maximal corners. All the information for doing that is given in section 5 and 6 of \cite{Yang2}. Thus we need to show that those are all maximal corners. We do this exemplary for case IV. Yang showed that $x=\sum_{i=1}^r \lambda_i x_i\in P_{\Gamma}$ if and only if
\begin{align}\label{Case IV verification}
\lambda_1 - \lambda_2\geq 0, \ldots, \lambda_{r-1}-\lambda_r, \lambda_r\geq 0\hspace{.3cm} \lambda_1\leq \frac{2}{(\psi,\psi)}, \hspace{.3cm} \sum_{i=1}^r \lambda_i\leq \frac{r}{(\psi,\psi)}.
\end{align} 
We set $D=\lbrace \lambda_1,\ldots,\lambda_r\in [0,\frac{2}{(\psi,\psi)}]\:|\: \sum_{i=1}^r \lambda_i\leq \frac{r}{(\psi,\psi)}\rbrace$. In the proof of Lemma 6.1. in \cite{Yang2} it is shown that for every $(\mu_1,\ldots,\mu_r)\in D$ with $\sum_{i=1}^r \mu_i= \frac{r}{(\psi,\psi)}$ the intersection $\lbrace \mu_i\:|\: 1\leq i \leq r \rbrace\cap (0,\frac{2}{(\psi,\psi)})$ has at most one element. If we take into account that by \eqref{Case IV verification} the elements $\lambda_i$ should be arranged in decreasing order and that $e_{\frac{r}{2}}=\sum_{i=1}^{\frac{r}{2}} \frac{2}{(\psi,\psi)} x_i$, $\frac{1}{2}(e_{\frac{r-1}{2}}+e_{\frac{r+1}{2}})=\sum_{i=1}^{\frac{r-1}{2}} \frac{2}{(\psi,\psi)} x_i + \frac{1}{(\psi,\psi)} x_{\frac{r+1}{2}}$, the uniqueness of the maximal corner follows.

\subsection*{Crittendens theorem}

Here we want to derive Theorem \ref{kernelexpcor} from Theorem 3 in \cite{Cr}.\\
Let $e^{\cdot}$ be the exponential map in $G$ and $P:=e^{\p}$. Then $P$ is a symmetric space and there is a diffeomorphism $\mu:G\slash K \to P$ defined by $\mu(gK):=g(\sigma(g))^{-1}$, where $\sigma$ is the involution on $G$ (see \cite{Cr} p. 321). It is well known that $(d\mu)_{eK}:T_{eK}G\slash K \cong \p\to\p$ is of the form $(d\mu)_{eK}=2\cdot Id$. Furthermore the following diagram is commutative (see \cite{Cr} p. 321)

\begin{center}
\begin{tikzpicture}
  \matrix (m) [matrix of math nodes,row sep=3em,column sep=4em,minimum width=2em]
  {
     T_{eK} G\slash K & \p \\
     G\slash K & P \\};
  \path[-stealth]
    (m-1-1) edge node [left] {$\exp_{eK}$} (m-2-1)
            edge  node [above] {$d\mu$} (m-1-2)
    (m-2-1.east|-m-2-2) edge node [below] {$\mu$}
             (m-2-2)
    (m-1-2) edge node [right] {$e^{\cdot}$} (m-2-2)
            (m-2-1);
\end{tikzpicture}
\end{center}

In the paper \cite{Cr} roots are defined as those elements $\bar{\alpha}\in\hp^{*}$, where a non-zero, 2-dimensional invariant subspace $V_{\bar{\alpha}}$ exists, such that for a suitable basis 
\begin{align*}
\ad_{\g} (\hp)_{|V_{\bar{\alpha}}}=\begin{pmatrix}
0 & 2\pi\bar{\alpha}\\ -2\pi\bar{\alpha}&0
\end{pmatrix}.
\end{align*}
We denote by $\bar{\Sigma}$ the set of these roots. Let
$\mathbf{V}_{\bar{\alpha}}:=\sum_{\bar{\beta}= \pm\bar{\alpha}} V_{\bar{\beta}}$, while for $\bar{\beta}$ there might exist more than one $V_{\bar{\beta}}$. 

\bl \label{root-space-relation}\mbox{}\\
Notations as before. Then $\textbf{V}_{\bar{\alpha}}= \g(2\pi\bar{\alpha})$ and $\bar{\Sigma}=\frac{\Sigma}{2\pi}$.
\el

\begin{proof}
Since
\begin{align*}
\ad_{\g} (\hp)_{|V_{\bar{\alpha}}}^2=\begin{pmatrix}
-4\pi^2\bar{\alpha}^2&0\\ 0&-4\pi^2\bar{\alpha}^2
\end{pmatrix},\hspace{.5cm} \ad_{\g} (\hp)_{|V_{-\bar{\alpha}}}^2=\begin{pmatrix}
-4\pi^2\bar{\alpha}^2&0\\ 0&-4\pi^2\bar{\alpha}^2
\end{pmatrix}
\end{align*}
it follows that $\mathbf{V}_{\bar{\alpha}}\subseteq \g(2\pi\bar{\alpha})$. The root space decomposition of $\g$ and the decomposition of $\g$ into eigenspaces of $\bar{\alpha}\in\bar{\Sigma}$ (see \cite{Cr} p. 322) are of the form 
\begin{align*}
\g=\g(0)\oplus\bigoplus_{\alpha\in\Sigma^{+}}\g(\alpha)\quad \text{and}\quad \g=\g(0)\oplus\bigoplus_{\bar{\alpha}\in\bar{\Sigma^{+}}}\mathbf{V}_{\bar{\alpha}},
\end{align*}

respectively. By dimensional reasons it has to be $\textbf{V}_{\bar{\alpha}}= \g(2\pi\bar{\alpha})$. This gives $2\pi\bar{\Sigma}=\Sigma$ or $\bar{\Sigma}=\frac{\Sigma}{2\pi}$.
\end{proof}

Since $\mathbf{V}_{\bar{\alpha}}= \g(2\pi\bar{\alpha})$, it follows from the definition of $\p(\alpha)$ that $\mathbf{V}_{\bar{\alpha}}\cap\p = \g(2\pi\bar{\alpha})\cap\p=\p(2\pi\bar{\alpha}) =\p(\alpha)$ for $\alpha\in \Sigma$ with $\alpha=2\pi\bar{\alpha}$. Furthermore let $n:T_0 \p\to \p$ denote the natural identification. We cite the theorem which we want to apply to our situation.

\bt\label{kernelexp1}(see \cite{Cr} p. 325)\\
Let $H\in\hp$. Then the Euclidean parallel translate of 
\begin{align*}
n^{-1}(\bigoplus\limits_{\bar{\alpha}(H)\equiv 0\mod 1\: ;\:\bar{\alpha}(H)\neq 0\: ;\: \bar{\alpha}\in\bar{\Sigma}^{+}}\mathbf{V}_{\bar{\alpha}} \cap \p\hspace{.2cm}) 
\end{align*}  
to $H$ constitutes the kernel of the differential $e^{\cdot}$ at $H$.
\et

Now we are able to derive Theorem \ref{kernelexpcor}.

\begin{proof}
Since $\mu$ is a diffeomorphism, we have $\ker((d\exp_{eK})_H)=\ker((d\mu\circ d\exp_{eK})_H)$. With the commutativity of the diagram above it follows 
\begin{align*}
\ker((d\exp_{eK})_H)=\ker(d( e^{\cdot}\circ d\mu))_H) = \ker((de^{2\cdot})_H).
\end{align*}
If we apply Theorem~\ref{kernelexp1}, it follows that those spaces $\textbf{V}_{\bar{\alpha}} \cap \p$ constitute to $\ker((de^{2\cdot})_H)$, where $\bar{\alpha}\in\bar{\Sigma}^{+}$ with $\bar{\alpha}(2H)\equiv 0\mod 1$ and $\bar{\alpha}(2H)\neq 0$. In view of Lemma \ref{root-space-relation}, we can write $\bar{\alpha}=\frac{\alpha}{2\pi}$ with $\alpha\in\Sigma$. Hence 
\[\bar{\alpha}(2H)=\frac{2\alpha(H)}{2\pi}\equiv 0\mod 1 \:\Longleftrightarrow\: \alpha(H)=0\mod \pi. \]
In the same way $\bar{\alpha}(2H)\neq 0$ is equivalent to $\alpha (H)\neq 0$. As $\textbf{V}_{\bar{\alpha}} \cap \p=\p(2\pi\bar{\alpha})=\p(\alpha)$, we get
\begin{align*}
\ker((d\exp_{eK})_H)= \ker((de^{2\cdot})_H)=& n^{-1}(\bigoplus\limits_{\bar{\alpha}(2H)\equiv 0\mod 1\: ;\:\bar{\alpha}(2H)\neq 0\: ;\: \bar{\alpha}\in\bar{\Sigma}^{+}}\textbf{V}_{\bar{\alpha}} \cap \p\hspace{.2cm})\\
=& n^{-1}(\bigoplus\limits_{\alpha(H)\equiv 0\mod \pi\: ;\:\alpha(H)\neq 0\: ;\: \alpha\in\Sigma^{+}}\p(\alpha)\hspace{.2cm}), 
\end{align*}
for $\alpha\in\Sigma^{+}$ with $2\pi\bar{\alpha}=\alpha$.
\end{proof}

For simplicity the natural identification $n$ was omitted in the sections above.

\vskip .5cm

Author's address:\\

Jonas Beyrer\\
PhD Student\\
Mathematik Institut, Univerist\"at Z\"urich,\\
Winterthurerstrasse 190, 8057 Z\"urich, Switzerland\\
E-mail: jonas.beyrer@math.uzh.ch


\begin{thebibliography}{9}

\bibitem{Bou2}{}
	N. Bourbaki:
	Lie Groups and Lie Algebras Chapters 4 - 6, Springer, Berlin-Heidelberg, 1st English ed., 1st softcover printing, 2008.	


\bibitem{Cr}{}
	R. J. Crittenden:
	\emph{Minimum and conjugate points in symmetric spaces}, Canad. J. Math. \textbf{14} (1962), 320-328.

\bibitem{DL}{}
	S. Deng, X. Liu: \emph{The antipodal sets of compact symmetric spaces}, Balk. J. Geom. Appl. \textbf{19} No.1 (2014), 73-79.

\bibitem{Hel}{}
  S. Helgason:
  Differential Geometry, Lie groups, and Symmetric Spaces,
  Graduate studies in Mathematics \textbf{34}, Amer. Math. Soc., 2000.
	
\bibitem{Kon}{}
	K. Kondo:
	\emph{Local orbits of \textit{S}-representations of symmetric $\R$-spaces},
	Tokyo J. Math \textbf{26} No. 1 (2003), 67-81. 	

\bibitem{Lo}{}
  O. Loos:
  Symmetric Spaces II: Compact Spaces and Classification,
  W.A. Benjamin, INC, 1969. 
  
\bibitem{Ti}{}
	J. A. Tirao: \emph{Antipodal manifolds in compact symmetric spaces of rank one},
	Proc. Amer. Math. Soc. \textbf{72} No. 1 (1978), 143–149. 
	
\bibitem{Yang}{}
	L. Yang:
	\emph{Injectivity radius and Cartan polyhedron for simply connected symmetric spaces}, Chinese Ann. Math. B \textbf{28} Issue 6 (2007), 685-700.
	
\bibitem{Yang2}{}
	L. Yang: \emph{Injectivity radius for non-simply connected symmetric spaces via Cartan polyhedron}, Osaka J. Math. \textbf{45} No. 2 (2008), 511-540.
\end{thebibliography}
\end{document}